\documentclass[12pt,twoside]{amsart}\usepackage{amssymb}
\usepackage{enumerate}

\newtheorem{thmspec}{\relax}

\newtheorem{theorem}{Theorem}[section]
\newtheorem{thm}[theorem]{Theorem}
\newtheorem{lem}[theorem]{Lemma}

\newtheorem{prop}[theorem]{Proposition}
\newtheorem{defi}[theorem]{Definition}
\newtheorem{rem}[theorem]{Remark}
\theoremstyle{definition}

\theoremstyle{remark}

\numberwithin{equation}{section}
\def \Bbb{\mathbb}

\def\onto{{\kern3pt\to\kern-8pt\to\kern3pt}}

\def\<{\langle}
\def\>{\rangle}
\def\|{{\ |\ }}

\def\onto{\twoheadrightarrow}

\def\-{\underline}

\def\N{\Bbb N}

\def\R{\Bbb R}
\def\C{\Bbb C}
\def\P{\Bbb P}

\def\X{\Bbb X}

\def\<{\langle}
\def\>{\rangle}

\catcode`\@=11
\def\serieslogo@{\relax}
\def\@setcopyright{\relax}
\catcode`\@=12

\title[ A general version of the Hartogs extension theorem]
{   A  general version of the Hartogs extension\\
 theorem for  separately holomorphic mappings\\ between   complex analytic spaces
}

\begin{document}

\author{Vi\^et-Anh  Nguy\^en}
\address{Vi\^et-Anh  Nguy\^en\\
 Max-Planck-Institut f\"{u}r Mathematik   \\
Vivatsgasse 7\\
 D--53111  Bonn\\
 Germany}
\email{vietanh@mpim-bonn.mpg.de}

\subjclass[2000]{Primary 32D15, 32D10}
\date{}

\keywords{  Cross theorem,  holomorphic extension,   Poletsky theory of discs,
Rosay's Theorem on holomorphic discs.}

\begin{abstract}
Using recent development in   Poletsky theory of discs, we prove the following result:
  Let $X,$ $Y$ be two complex manifolds, let  $Z$ be
  a complex analytic  space which possesses the Hartogs extension property, let
  $A$ (resp. $B$) be a  non locally pluripolar subset of  $X$ (resp.
  $Y$). We   show that every separately
  holomorphic
  mapping  $f:\  W:=(A\times Y) \cup (X\times B)\longrightarrow Z$ extends to a    holomorphic mapping
  $\hat{f}$ on  $\widehat{W}:=\left\lbrace(z,w)\in X\times Y:\ \widetilde{\omega}(z,A,X)+\widetilde{\omega}(w,B,Y)<1  \right\rbrace$
  such that  $\hat{f}=f$ on $W\cap \widehat{W},$
   where $\widetilde{\omega}(\cdot,A,X)$ (resp. $\widetilde{\omega}(\cdot,B,Y))$ is the plurisubharmonic measure of
   $A$ (resp. $B$) relative to $X$ (resp. $Y$).   
   Generalizations of this result for an $N$-fold cross   are also given.
\end{abstract}
\maketitle

\section{ Introduction}
The main purpose of this article  is to give a general version  of the
well-known Hartogs extension theorem for separately holomorphic functions
(see \cite{ha}). This theorem has been a source of inspiration for numerous
research works in Complex Analysis for many years. It has developed into
the beautiful and very active theory of     separately analytic
 mappings.   Nowadays, one finds    a close connection  between this theory
 and
many other fields in Mathematics such as (Pluri)potential Theory, Partial Differential Equations
and Theoretic Physics  etc.
The recent survey
articles by  Nguy\^en Thanh V\^an \cite{ng} and Peter  Pflug \cite{pf} not only retrace
the historic development, but also give some insights into the
new research trends in this subject. Here
we  recall    briefly  the main steps  in  developing the theory of separately
holomorphic  mappings. 

\smallskip

Very longtime after the ground-breaking work of Hartogs, the subject was re-birthed,  around the years 50--60s,
thanks to the Japanese school (see \cite{sh}, \cite{te} and the references
therein). However, an important impetus  was only made by Siciak in  the works \cite{si1,si2}, where he established some significant  generalizations  of
the Hartogs
extension theorem.  According to  Siciak's general  formulation of this
theorem,
the problem is  to determine the envelope of holomorphy for  separately
holomorphic functions defined on some {\it cross sets.} The theorems
obtained under this formulation are often called {\it cross theorems.}
 Using  the so-called {\it relative
extremal function,}  Siciak completed the problem for  the case where
the cross set consists of a product of domains in $\C.$

\smallskip

The next deep steps were initiated by Zahariuta in 1976 (see \cite{za}) when he started to use
the method of common bases of Hilbert spaces. This original approach permits
him to obtain new  cross theorems for some cases where the cross consists
of a product of Stein manifolds. As a consequence, he was able to generalize
the result of Siciak in higher dimensions.

\smallskip

Later, Nguy\^en  Thanh V\^an and Zeriahi (see \cite{nz1,nz2,nz3}) developed the method of doubly orthogonal bases of
Bergman type in order to generalize the result of
Zahariuta. This  is a significantly simpler and    more constructive version
of Zahariuta's original  method.  Nguy\^en  Thanh V\^an and Zeriahi have recently achieved an elegant
 improvement  of their method  (see \cite{ng2}, \cite{ze}).

\smallskip

Using the method of Siciak, Shiffman  (see \cite{sh1}) was the first to generalize some
results of Siciak to separately holomorphic mappings with values in a complex analytic space.



\smallskip

The most general result  to date is contained in a recent work by  Alehyane and Zeriahi
(see Theorem 2.2.4 in  \cite{az}). Namely, they are able to define the envelope of holomorphy of any cross of a
product of subdomains of Stein manifolds  in terms of the plurisubharmonic measure.

\smallskip

In this work we generalize, in some sense, the result of  Alehyane--Zeriahi  to any cross of a product of arbitrary
complex manifolds. The main ingredient in our approach is  Poletsky theory of   discs developed in \cite{po1, po2},
    Rosay's Theorem on holomorphic discs (see \cite{ro}) and  Alehyane--Zeriahi Theorem
(see \cite{az}). Another important technique is to use level sets of
the plurisubharmonic measure. This technique was originally introduced in a recent
joint-work of   Pflug and the author (see \cite{pn1}). 
However, it appears to be very successful in solving many problems arising  from the
theory of separately holomorphic and meromorphic mappings.

\smallskip

This paper is organized as follows.

\smallskip

 In Section 2,  after   introducing some terminology and notation, we
 recall  Alehyane--Zeriahi  Theorem and state our main result.

\smallskip

The tools which are needed for the proof of the main result are developed
in Sections 3.

\smallskip

   The proof of the    main result for  the   case of an $2$-fold cross   is divided into three
   parts, which correspond to Sections 4, 5, and 6.


\smallskip

   The general case  is treated in Section 7.
  Finally, we conclude the article with  some   remarks and open questions.

 \smallskip

 The theory of separately holomorphic and meromorphic mappings has
  received  much attention  in the past few years.
 We only mention here some directions of the current research.
 Separate analyticity in infinite dimension is growing quite rapidly
since the work of Noverraz \cite{no}. Many results in this direction are obtained   by
  Nguy\^en Van  Khu\^e,   Nguy\^en  Thanh V\^an  and their
co-workers (see the discussion in \cite{kt} and \cite{ng}).
On the other hand, the recent development also focuses on cross theorems with pluripolar
singularities and  boundary cross theorems. For the latest results as well as
a comprehensive introduction to the latter two directions,
the reader may consult  some works of Jarnicki and Pflug
in
\cite{jp2,jp3,jp4}
and   recent articles of  Pflug and the author  (see \cite{pn1,pn2}).

\smallskip

\indent{\it{\bf Acknowledgment.}} The paper was written while  the  author was visiting the
Carl von  Ossietzky Universit\"{a}t Oldenburg
 being  supported by The Alexander von Humboldt
Foundation. He wishes to express his gratitude to these organizations. He
also would like to thank Professor Peter Pflug for very stimulating
discussions, and the referee for many valuable remarks.

\section{Preliminaries and statement of the main result}
In order to recall  the classical  cross theorem and to
 state the main result, we need to introduce some notation and terminology.
 In fact, we keep the main notation from the   works in \cite{jp1},  \cite{sa}.
 \subsection{Local pluripolarity and plurisubharmonic measure, cross and  separate holomorphicity.}
 In the sequel, all complex manifolds are supposed to be  of  finite local dimension (i.e. the dimension
 of any   connected component of the manifold is finite), and all
  complex analytic spaces considered in this work are supposed to be reduced,
  irreducible and of finite dimension.

  Let $\mathcal{M}$ be a complex manifold  and  let $A$ be a  subset
   of $ \mathcal{M}.$  Put
\begin{equation*}
h_{A,\mathcal{M}}:=\sup\left\lbrace u:\ u\in\mathcal{PSH}(\mathcal{M}),\ u\leq 1 \ \text{on}\
\mathcal{M}, \ u\leq 0\ \text{on}\ A \right\rbrace,
\end{equation*}
where $\mathcal{PSH}(\mathcal{M})$ denotes the cone of all plurisubharmonic
functions on $\mathcal{M}.$

$A$ is said to be {\it pluripolar} in $\mathcal{M}$ if there is  $u\in
\mathcal{PSH}(\mathcal{M})$ such that $u$ is not identically $-\infty$ on every
connected component of $\mathcal{M}$ and $A\subset \left\lbrace z\in \mathcal{M}:\ u(z)=-\infty\right\rbrace.$
$A$ is said to be {\it locally  pluripolar} in $\mathcal{M}$ if  for any
$z\in A,$ there is an open  neighborhood $V$ of $z$ such that $A\cap V$ is
pluripolar in $V.$ $A$ is said  to be {\it nonpluripolar} (resp. {\it non locally  pluripolar}) if it is not
pluripolar (resp. not locally pluripolar).
According to a classical result of Josefson and Bedford  (see \cite{jo},
\cite{be}), if $\mathcal{M}$ is a Riemann domain over a Stein manifold,
then   $A\subset \mathcal{M}$ is   locally  pluripolar if and only if it is
pluripolar.

In the sequel, for a function $h:\ \mathcal{M}\longrightarrow \R,$ its
{\it  upper semicontinuous regularization $h^{\ast} :\ \mathcal{M}\longrightarrow \R $} is defined by
\begin{equation*}
h^{\ast}(z):=\limsup\limits_{w\to z}h(w),\qquad  z\in\mathcal{M}.
\end{equation*}

 Next, we say that a  set  $A\subset \mathcal{M}$ is {\it
locally pluriregular at} a point $a\in \overline{A}$   if $h^{\ast}_{A\cap U,U}(a)=0$
  for  all open neighborhoods $U$ of $a.$
Moreover, $ A$ is said to be {\it locally pluriregular } if it is locally
pluriregular at every point $a\in A.$  We denote by $A^{\ast}=A^{\ast}_{\mathcal{M}}$ the set of all
points  $a\in \overline{A} $ at which $A$ is locally pluriregular.
If $A$ is non locally pluripolar, then  a classical result of Bedford and Taylor
(see \cite{be,bt})
says that $A^{\ast}$ is non locally pluripolar and  $A\setminus A^{\ast}$ is locally  pluripolar.
Moreover, $A^{\ast}$ is locally  of type $\mathcal{G}_{\delta}$
(i.e. for every $a\in A^{\ast},$ there is an open  neighborhood $U$ of $a$
such that $A^{\ast}\cap U$ is  a countable intersection of open sets), and  $A^{\ast}$ is locally
pluriregular (i.e. $(A^{\ast})^{\ast}=A^{\ast}$).

 The {\it  plurisubharmonic measure} of $A$ {\it relative to $\mathcal{M}$} is
 the function\footnote{ The  notation $\omega(\cdot,A,\mathcal{M})$
is historically reserved for the
  {\it relative extremal function}.  The latter function is  defined by   $$\omega(z,A,\mathcal{M}):= h^{\ast}_{A,
\Omega}(z),\qquad  z\in \mathcal{M}.$$  An example  in  \cite{ah} shows  that     in general,  $\omega(\cdot,A,\mathcal{M})\not=\widetilde{\omega}(\cdot,A,\mathcal{M}).$} defined by
\begin{equation*}
\widetilde{\omega}(z,A,\mathcal{M}):=  h^{\ast}_{A^{\ast},
\Omega}(z),\qquad  z\in \mathcal{M}.
\end{equation*}
Observe that  $\widetilde{\omega}(\cdot,A,\mathcal{M})\in\mathcal{PSH}(\mathcal{M})$ and
$0\leq \widetilde{\omega}(z,A,\mathcal{M})\leq 1,\ z\in\mathcal{M}.$

We shall show in Proposition \ref{prop3.6} below  that if $\mathcal{M}$ is
a subdomain of a Stein manifold, then
the above  definition coincides  with the one given by formula (2.1.2) in  Alehyane--Zeriahi's article
 \cite{az}.
 For a good background of the pluripotential
theory, see the books  \cite{jp1} or  \cite{kl}.

 Let $N\in\N,\ N\geq 2,$ and let $\varnothing\not = A_j\subset
D_j,$ where $D_j$ is a  complex manifold, $j=1,\ldots,N.$ We define
an {\it $N$-fold cross}
\begin{eqnarray*}
X &:=&\X(A_1,\ldots,A_N; D_1,\ldots,D_N)\\
&:=& \bigcup_{j=1}^N  A_1\times\cdots\times A_{j-1}\times D_j\times
A_{j+1}\times\cdots A_N .
\end{eqnarray*}
Following a terminology of Alehyane--Zeriahi \cite{az}, we define the {\it regular
part} $X^{\ast}$ of $X$ as follows
\begin{eqnarray*}
X^{\ast} &=&\X^{\ast}(A_1,\ldots,A_N; D_1,\ldots,D_N):=
\X(A^{\ast}_1,\ldots,A^{\ast}_N; D_1,\ldots,D_N)\\
&=& \bigcup_{j=1}^N  A^{\ast}_1\times\cdots\times A^{\ast}_{j-1}\times D_j\times
A^{\ast}_{j+1}\times\cdots A^{\ast}_N .
\end{eqnarray*}
Moreover, put
\begin{equation*}
\omega(z):=\sum\limits_{j=1}^N \widetilde{\omega}(z_j,A_j,D_j),\qquad
z=\left(z_1,\ldots,z_N\right)\in D_1\times
\cdots\times D_N.
\end{equation*}
For an $N$-fold cross $X:=\X(A_1,\ldots,A_N;D_1,\ldots,D_N)$ let
\begin{equation*}
\widehat{X}=\widehat{\X}(A_1,\ldots,A_N;D_1,\ldots,D_N) :=\left\lbrace \left(z_1,\ldots,z_N\right )\in D_1\times
\cdots\times D_N:\  \omega(z)<1
\right\rbrace.
\end{equation*}
It is not difficult to see that $X^{\ast}\subset \widehat{X}.$

Let $Z$ be a complex analytic space.
We say that a mapping $f:X\longrightarrow Z$ is {\it separately holomorphic}
  and write $f\in\mathcal{O}_s(X,Z),$   if
 for any $j\in\lbrace
1,\ldots,N\rbrace$ and $(a^{'},a^{''})\in  (A_1\times\cdots\times
 A_{j-1})\times(A_{j+1}\times \cdots \times A_N) $
 the restricted mapping $f(a^{'},\cdot,a^{''})|_{D_j}$ is holomorphic  on $D_j.$

Throughout the paper, for a  function  $f:\ M\longrightarrow\C,$ let
   $\vert f\vert_M$ denote $\sup_M \vert f\vert.$
Finally, for a complex manifold $\mathcal{M}$ and a complex analytic space $Z,$ let
$\mathcal{O}(\mathcal{M},Z)$ denote the set of all holomorphic mappings from $\mathcal{M}$
to $Z.$  

\subsection{Hartogs extension property.} We recall here  the following notion
introduced by Shiffman \cite{sh1}.
 Let $p\geq 2$ be an integer. For $0<r<1,$  the {\it Hartogs
 figure} in dimension $p,$ denoted by $H_p(r),$ is given by
 \begin{equation*}
H_p(r):=\left\lbrace (z^{'},z_p)\in E^p: \  \Vert z^{'}\Vert<r \ \ \text{or}\ \ \vert z_p\vert >1-r \right\rbrace,
 \end{equation*}
where $E$ is the unit disc of $\C$ and $z^{'}=(z_1,\ldots,z_{p-1}),$
$\Vert z^{'}\Vert:=\max\limits_{1\leq j\leq p-1} \vert z_j\vert.$

\renewcommand{\thethmspec}{Definition 1}
  \begin{thmspec}
A complex analytic space $Z$ is said to possess the Hartogs extension property
  in dimension $p$ if $Z$ has a countable bases of open subsets, and every  mapping $f\in\mathcal{O}(H_p(r) , Z)$ extends to a  mapping
  $\hat{f}\in\mathcal{O}(E^p , Z).$  Moreover, $Z$ is said to possess the Hartogs extension property
 if it does in any dimension $p\geq 2.$
\end{thmspec}

It is a classical result of Ivashkovich (see \cite{iv1}) that if $Z$ possesses
the Hartogs extension property  in dimension 2, then it does in all dimensions $p\geq 2.$
 Some typical  examples of complex analytic  spaces  possessing the Hartogs extension property are  the
complex Lie groups (see \cite{asy}), the taut spaces (see \cite{wu}), the Hermitian manifold with
negative holomorphic
sectional  curvature (see \cite{sh1}), the holomorphically convex K\"{a}hler  manifold without rational curves
(see \cite{iv1}).

Here we  mention  an important characterization due to Shiffman (see \cite{sh1}).
\renewcommand{\thethmspec}{Theorem 2}
  \begin{thmspec}
A complex analytic space $Z$   possesses the Hartogs extension property
    if and only if  for every domain $D$ of any Stein manifold $\mathcal{M},$ every mapping
    $f\in\mathcal{O}(D, Z)$ extends to a  mapping $\hat{f}\in
    \mathcal{O}(\widehat{D},Z) ,$   where $\widehat{D}$ is the envelope of holomorphy
of $D.$
\end{thmspec}
\subsection{Motivations for our work}
We are now able to formulate what we will quote in the sequel as the {\it classical cross
theorem}.
\renewcommand{\thethmspec}{Theorem 3}
  \begin{thmspec} (Alehyane--Zeriahi \cite[Theorem 2.2.4]{az})
  Let $X_j$ be a  Stein manifold, let  $D_j\subset X_j$ be a    domain
  and  $ A_j\subset
D_j$ a nonpluripolar   subset,
  $j=1,\ldots,N.$   Let $Z$ be a complex  analytic space possessing the
  Hartogs
  extension property.
Then for any mapping $f\in\mathcal{O}_s(X,Z),$ there is a unique mapping
$\hat{f}\in\mathcal{O}(\widehat{X},Z)$
such that $\hat{f}=f$ on $X\cap \widehat{X}.$
\end{thmspec}

The following example  given  by Alehyane--Zeriahi (see \cite{az}) shows that the
hypothesis on $Z$ is necessary. Consider the mapping $f:\
\C^2\longrightarrow \P^1$ given by
\begin{equation*}
  f(z,w):=
\begin{cases}
[(z+w)^2: (z-w)^2],
  & (z,w)\not=(0,0) ,\\
 [1:1], &  (z,w)=(0,0).
\end{cases}
\end{equation*}
Then $f\in\mathcal{O}_s\Big(\X(\C,\C;\C,\C),\P^1\Big),$ but $f$ is not
continuous at $(0,0).$

The question naturally arises whether Theorem 3 is still true if
 $D_j$ is not necessarily a subdomain of a Stein manifold, $j=1,\ldots,N.$

\subsection{Statement of the  main results and outline of the proofs}
We are now ready to state the   main results.

\renewcommand{\thethmspec}{Theorem A}
  \begin{thmspec}
Let  $D_j $ be a complex manifold
  and  $ A_j\subset
D_j$ a non locally pluripolar   subset,
  $j=1,\ldots,N.$ Let $Z$ be a complex analytic space possessing the
  Hartogs
  extension property.
Then for any mapping $f\in\mathcal{O}_s(X,Z),$ there is a unique mapping
$\hat{f}\in\mathcal{O}(\widehat{X},Z)$
such that $\hat{f}=f$ on $X\cap \widehat{X}.$
If, moreover,  $Z=\C$ and $\vert f\vert_X<\infty,$ then
\begin{equation*}
 \vert \hat{f}(z)\vert\leq \vert f\vert_A^{1-\omega(z)} \vert
 f\vert_X^{\omega(z)},\qquad z\in\widehat{X}.
\end{equation*}
\end{thmspec}

In virtue of a theorem of  Josefson and Bedford (see Subsection 2.1
above),  the classical cross theorem is an immediate consequence of Theorem A.

Theorem A has an important corollary. Before stating this, we need to
introduce a terminology. A complex manifold $\mathcal{M}$ is said to be
 a {\it Liouville manifold} if $\mathcal{PSH}(\mathcal{M})$ does not
 contain any non-constant bounded above functions. We see clearly that the class
 of Liouville manifolds contains the class of connected compact manifolds.

\renewcommand{\thethmspec}{Corollary B}
  \begin{thmspec}
Let  $D_j $ be a complex manifold
  and  $ A_j\subset
D_j$ a non locally pluripolar   subset,
  $j=1,\ldots,N.$ Let $Z$ be a complex analytic  space possessing the  Hartogs extension property.
   Suppose in addition that  $D_j $ is  a Liouville manifold,
  $j=2,\ldots,N.$
Then for any mapping $f\in\mathcal{O}_s(X,Z),$ there is a unique mapping
$\hat{f}\in\mathcal{O}(D_1\times\cdots\times D_N,Z)$
such that $\hat{f}=f$ on $X .$
\end{thmspec}

Corollary B follows immediately from Theorem A since
$\widetilde{\omega}(\cdot,A_j,D_j)\equiv 0,$  $j=2,\ldots,N.$

We give below some ideas  of the proof of   Theorem A.

\smallskip

 Our method consists in two steps.
In the first step, we   investigate the special case where each $A_j$ is  an open set, $j=1,\ldots,N.$
In the second one, we treat the general case.

 In order to carry out the first step,
 we apply Poletsky theory of   discs  and Rosay's Theorem on holomorphic discs (see Theorem \ref{Rosaythm} below).
 Consequently, we may construct an extension mapping $\hat{f}$ on $\widehat{X}.$
 To prove that $\hat{f}$ is holomorphic, we appeal to the classical cross theorem (Theorem 3).

In the second step we   reduce the general situation to the above
special case. The key technique is to
use {\it level sets} of the plurisubharmonic measure. More precisely,
 we exhaust each $D_j$ by the  level sets of the plurisubharmonic measure
$\widetilde{\omega}(\cdot,A_j,D_j),$ i.e. by
$D_{j,\delta}:=\left\lbrace z_j\in D_j:\ \widetilde{\omega}(z_j,A_j,D_j)<1-\delta \right\rbrace$
$(0<\delta<1).$  We replace in the same way  the set  $A_j$ by an open set
$A_{j,\delta}$ such that
$\widetilde{\omega}(\cdot,A_{j,\delta},D_{j,\delta} )$ behaves, in some sense, like
$\widetilde{\omega}(\cdot,A_j,D_j)$ as $\delta\to 0^{+}.$
  Applying  Theorem  3 locally and  making an intensive use
  of  Theorem \ref{Rosaythm}, we can propagate the separate holomorphicity
of $f$ to a   mapping $\tilde{f}_{\delta}$  defined on the cross
 $X_{\delta}:=\X\left(A_{1,\delta},\ldots,A_{N,\delta};D_{1,\delta},\ldots,
 D_{N,\delta}\right).$ Consequently, the first step applies and one
 obtains a mapping $\hat{f}_{\delta}\in\mathcal{O}\left(\widehat{X_{\delta}},Z\right). $
 Gluing the family $\left(\hat{f}_{\delta}\right)_{0<\delta<1},$ we obtain the desired
 extension mapping $\hat{f}.$

\section{Preparatory results}
We recall here the auxiliary  results and some  background of the
pluripotential theory needed for the proof of Theorem A.
\subsection{Poletsky theory of  discs and Rosay's Theorem on
 holomorphic  discs}
  Let $E$ denote as usual the unit disc in $\C.$ %
For a complex
manifold $\mathcal{M},$ let $\mathcal{O}(\overline{E},\mathcal{M})$ denote
the set of all holomorphic mappings $\phi:\ E\longrightarrow \mathcal{M}$ which
extend holomorphically  to   a neighborhood of  $\overline{E}.$
Such a mapping $\phi$ is called a {\it holomorphic disc} on $\mathcal{M}.$ Moreover, for
a subset $A$ of $\mathcal{M},$ let
\begin{equation*}
 1_{  A}(z):=
\begin{cases}
1,
  &z\in   A,\\
 0, & z\in \mathcal{M}\setminus A.
\end{cases}
\end{equation*}

In the work \cite{ro}  Rosay proved the following remarkable result.
\begin{thm}\label{Rosaythm}
Let $u$ be an upper semicontinuous function on a complex manifold
$\mathcal{M}.$ Then the Poisson functional of $u$  defined by
\begin{equation*}
\mathcal{P}[u](z):=\inf\left\lbrace\frac{1}{2\pi}\int\limits_{0}^{2\pi} u(\phi(e^{i\theta}))d\theta:  \
\phi\in   \mathcal{O}(\overline{E},\mathcal{M}), \ \phi(0)=z
\right\rbrace,
\end{equation*}
is plurisubharmonic on $\mathcal{M}.$
\end{thm}

Rosay's Theorem may be viewed as an important development in Poletsky
theory of   discs. Observe that special cases of Theorem
\ref{Rosaythm} have been considered by Poletsky (see \cite{po1,po2}),
L\'arusson--Sigurdsson (see \cite{ls}) and Edigarian (see \cite{ed}).

We also need the following result  (see \cite[Lemmas 1.1 and 1.2]{ro}).
\begin{lem}\label{Rosaylem}
Let $\mathcal{M}$ be a complex manifold and let $A$ be a nonempty open
subset of $\mathcal{M}.$   Then
for any $\epsilon>0$ and any $z_0\in \mathcal{M},$  there are an open neighborhood $U$ of
$z_0,$ an open subset $T$ of $\C,$    and a  family of holomorphic discs
 $(\phi_z)_{z\in U}\subset \mathcal{O}(\overline{E} ,\mathcal{M})$ with the
following properties:\\
(i) $\Phi\in\mathcal{O}(U\times E ,\mathcal{M}),$ where
$\Phi(z,t):=\phi_z(t),$  $(z,t)\in U\times E ;$\\
(ii) $\phi_{z}(0)=z,\qquad z\in U;$\\
(iii) $\phi_z(t)\in  A, \qquad t\in T\cap  \overline{E},\ z\in U;$\\
(iv)$  \frac{1}{2\pi}\int\limits_{0}^{2\pi}1_{ \partial E \setminus T } (e^{i\theta})d\theta
< \mathcal{P}[1_{\mathcal{M}\setminus A}](z_0 )+\epsilon.$
\end{lem}
\begin{proof} For any $\rho>0,$ let
 $E_{\rho}$ denote the disc  $\left\lbrace   t\in\C:\ \vert t\vert<\rho \right\rbrace.$ Fix an arbitrary point $z_0\in\mathcal{M}.$ We apply Theorem \ref{Rosaythm} to  the upper
semicontinuous function $1_{\mathcal{M}\setminus A}.$ Consequently, for any $\epsilon>0,$   one may
find an $r>1$ and a holomorphic mapping $\phi\in \mathcal{O}(E_r,\mathcal{M})$
such that
\begin{equation}\label{eq3.1.1}
\phi(0)=z_0\qquad \text{and}\qquad \frac{1}{2\pi}\int\limits_{0}^{2\pi}1_{ \mathcal{M}
\setminus A  } (\phi(e^{i\theta}))d\theta
<  \mathcal{P}[1_{\mathcal{M}\setminus A}](z_0)+\frac{\epsilon}{2}.
\end{equation}
  Consider the embedding $\tau:\
E_r\longrightarrow \C\times \mathcal{M}$  given by
$\tau (t):=(t,\phi(t)),$ $t\in  E_r.$ Then the image $\tau (  E_r)$ is
a Stein submanifold of $\C\times \mathcal{M}.$ Fix any $\tilde{r}$ such that
$1<\tilde{r}<r$ and let $d$ be the dimension of the connected component  of $\mathcal{M}$ containing $z_0.$
  By Lemma 1.1 in \cite{ro},
there is an injective holomorphic mapping $\tilde{\tau}:\
E_{\tilde{r}}^{d+1}\longrightarrow \C\times\mathcal{M}$ such that
$\tilde{\tau}(t,0)=\tau(t)=(t,\phi(t)), $  $\vert t\vert<\tilde{r}.$
Let $\Pi$ be the canonical projection from $\C\times\mathcal{M}$ onto
$\mathcal{M}.$ Then there are a sufficiently small neighborhood $U$ of
$z_0$ and a real number $\rho:$  $1<\rho<\tilde{r}$ such that, for every $z\in
U,$ the mapping $\phi_z:\  E_{\rho}\longrightarrow \mathcal{M}$ given
by
 \begin{equation}\label{eq3.1.2}
\phi_z(t):=\Pi\circ\tilde{\tau}\left((t,0)+\tilde{\tau}^{-1}(0,z)\right),\qquad   t\in
E_{\rho},
\end{equation}
 is holomorphic.

 Using the explicit formula (\ref{eq3.1.2}), assertion (i) follows.
 Moreover, $\phi_z(0)=\Pi(0,z)=z$ for $z\in U,$ which proves assertion
 (ii). In addition,
 \begin{equation}\label{eq3.1.3}
\phi_{z_0}(t)=(\Pi\circ\tilde{\tau})(t,0)=\phi(t).
 \end{equation}

In virtue of  (\ref{eq3.1.2}), observe that as $z$ approaches $z_0$ in $U,$  $\phi_z$
 converges uniformly to $\phi_{z_0}$ on $\overline{E}.$ Consequently, by shrinking $U$
 if necessary, we
 may find an open subset $T$ of the open set $ \left\lbrace t\in E_{\rho}:\ \phi_{z_0}(t)\in A\right\rbrace$ such that
 assertion (iii) is fulfilled and
\begin{equation*}
  \frac{1}{2\pi}\int\limits_{0}^{2\pi}1_{  \partial E\setminus T} (e^{i\theta})d\theta
<  \frac{1}{2\pi}\int\limits_{0}^{2\pi}1_{ \mathcal{M}
\setminus A  } (\phi_{z_0}(e^{i\theta}))d\theta +\frac{\epsilon}{2}.
\end{equation*}
This, combined with the estimate in (\ref{eq3.1.1}) and (\ref{eq3.1.3}), implies assertion (iv).
Hence, the proof of the  lemma is complete.
\end{proof}
\subsection{The plurisubharmonic measure and its level sets}
We begin this subsection with the following simple but very useful result.
\begin{lem}\label{lem3.3}
Let $T$ be an open subset of $\overline{E}.$   Then
\begin{equation*}
\widetilde{\omega}(0,T\cap E,E)\leq \frac{1}{2\pi}\int\limits_{0}^{2\pi}1_{  \partial E\setminus T}
(e^{i\theta})d\theta.
\end{equation*}
\end{lem}
\begin{proof}
Observe that, by definition,
\begin{equation*}
\widetilde{\omega}(t,T\cap E,E)\leq \omega_E(t,T\cap\partial E),\qquad t\in E,
\end{equation*}
where $ \omega_E(t,T\cap\partial E )$ is the harmonic measure for $E$ (see
\cite[p. 96]{ra}). Since
\begin{equation*}
\omega_E(0,T\cap \partial E ) = \frac{1}{2\pi}\int\limits_{0}^{2\pi}1_{  \partial E\setminus T}
(e^{i\theta})d\theta,
\end{equation*}
the desired conclusion follows from the above  estimate.
\end{proof}
\begin{prop}\label{prop3.4}
Let $\mathcal{M}$ be a complex manifold and   $A$   a nonempty open
subset of $\mathcal{M}.$    Then
   $ \widetilde{\omega}(z,A,\mathcal{M}) = \mathcal{P}[1_{\mathcal{M}\setminus A}](z),$ $z\in\mathcal{M}.$
\end{prop}
\begin{proof} First,  since $A$ is open, it is clear that $A^{\ast}=A.$
In addition, applying Theorem \ref{Rosaythm} to $1_{\mathcal{M}\setminus
A}$ and using the explicit formula of  $\mathcal{P}[1_{\mathcal{M}\setminus
A}],$ we see that $\mathcal{P}[1_{\mathcal{M}\setminus
A}]\in\mathcal{PSH}(\mathcal{M}),$ $\mathcal{P}[1_{\mathcal{M}\setminus
A}]\leq 1$ and $\mathcal{P}[1_{\mathcal{M}\setminus
A}](z)=0,$  $z\in A.$ Consequently,
\begin{equation*}
 \mathcal{P}[1_{\mathcal{M}\setminus A}](z)\leq
 \widetilde{\omega}(z,A,\mathcal{M}),\qquad z\in\mathcal{M}.
\end{equation*}
To see the opposite inequality, let $u\in\mathcal{PSH}(\mathcal{M})$ such
that $u\leq 1$ and $u(z)\leq 0,\  z\in A.$ For any point
$z_0\in\mathcal{M}$ and any $\epsilon>0,$ by Theorem \ref{Rosaythm},
   there is a  holomorphic  disc $\phi\in \mathcal{O}(\overline{E},\mathcal{M})$
such that
\begin{equation}\label{eq3.2.1}
\phi(0)=z_0\qquad \text{and}\qquad \frac{1}{2\pi}\int\limits_{0}^{2\pi}1_{ \mathcal{M}
\setminus A  } (\phi(e^{i\theta}))d\theta
<  \mathcal{P}[1_{\mathcal{M}\setminus A}](z_0)+\epsilon.
\end{equation}
Consequently, by setting $\phi^{-1}(A):=\left\lbrace t\in\overline{E}:\ \phi(t)\in A \right\rbrace,$ we
obtain
\begin{equation*}
u(z_0)=(u\circ\phi) (0) \leq \widetilde{\omega}(0, \phi^{-1}(  A),E)
\leq \frac{1}{2\pi}\int\limits_{0}^{2\pi}1_{ \mathcal{M}
\setminus A  } (\phi(e^{i\theta}))d\theta,
\end{equation*}
where the first estimate is trivial and the second one follows from Lemma
\ref{lem3.3}.  This, combined with (\ref{eq3.2.1}), implies that
$u(z_0)< \mathcal{P}[1_{\mathcal{M}\setminus A}](z_0)+\epsilon.$
Since $u,\ \epsilon$ and $z_0$ are arbitrarily chosen, we conclude that
$ \widetilde{\omega}(z,A,\mathcal{M})\leq\mathcal{P}[1_{\mathcal{M}\setminus A}](z),$ $
z\in\mathcal{M}.$ This completes the proof.
\end{proof}

\begin{prop}\label{prop3.5}
Let $\mathcal{M}$ be a complex manifold and   $A$  a non locally pluripolar
subset of $\mathcal{M}.$ For $0<\epsilon<1,$ define the "$\epsilon$-level
set of $\mathcal{M}$ relative to $A$" as follows
\begin{equation*}
 \mathcal{M}_{\epsilon,A}:=
\left\lbrace z\in\mathcal{M}:\ \widetilde{\omega}(z,A,\mathcal{M})<1-\epsilon \right\rbrace.
\end{equation*}
Then:\\
1) For every locally pluripolar subset $ P$  of $\mathcal{M},$
   $(A\cup P)^{\ast}=A^{\ast}$  and
$\widetilde{\omega}(\cdot,A\cup P,\mathcal{M})=
\widetilde{\omega}(\cdot,A ,\mathcal{M}).$   $(A^{\ast})^{\ast}=A^{\ast}.$  If, moreover, $A$ is open, then $A^{\ast}=A.$\\
2) Let $\mathcal{N}$ be an open subset of $\mathcal{M}$ and   $B\subset
A\cap \mathcal{N}.$   Then  $\widetilde{\omega}(z,A,\mathcal{M})\leq
\widetilde{\omega}(z,B,\mathcal{N}),$  $z\in\mathcal{N}.$\\
3) Let $\mathcal{N}$ be a connected component of $\mathcal{M},$ then
$\widetilde{\omega}(z,A\cap\mathcal{N},\mathcal{N})=\widetilde{\omega}(z,A,\mathcal{M}),$
$z\in\mathcal{N}.$\\
4)  $\widetilde{\omega}(z,A\cap
 A^{\ast},\mathcal{M}_{\epsilon,A})=\frac{\widetilde{\omega}(z,A,\mathcal{M})}{1-\epsilon},$
$z\in\mathcal{M}_{\epsilon,A}.$\\
5) Every connected component of $\mathcal{M}_{\epsilon,A}$ contains a
non locally pluripolar subset of $A\cap A^{\ast}.$ If, moreover, $A$ is open, then
every connected component of $\mathcal{M}_{\epsilon,A}$ contains a
nonempty open subset of $A$.
\end{prop}
\begin{proof}
Part 1) is an immediate consequence of the following identity (see Lemma
3.5.3 in \cite{jp1})
\begin{equation*}
h^{\ast}_{A\cup P, U}=h^{\ast}_{A,U},
\end{equation*}
where $U$ is a bounded open subset of $\C^n,$ $A$ and $ P$ are subsets of $U,$  and $P $ is
pluripolar.

Part 2) and Part 3) are trivial using the definition of the plurisubharmonic measure.

Now we turn to Part 4). Observe that for any $a\in   A^{\ast},$
\begin{equation}\label{eq3.2.2}
 \widetilde{\omega} (a,A,\mathcal{M})=  \widetilde{\omega}(a,A^{\ast},\mathcal{M})=0,
\end{equation}
where the first equality follows from  the definition of the
plurisubharmonic measure, and the second one from Part 2) and the
assumption that $a\in  A^{\ast}.$  Hence, $  A^{\ast}\subset
\mathcal{M}_{\epsilon,A}.$ In addition, we have clearly that
$\frac{\widetilde{\omega} (z,A,\mathcal{M})}{1-\epsilon}\leq 1,$
$z\in\mathcal{M}_{\epsilon,A}.$ This, combined with (\ref{eq3.2.2}),
implies that
\begin{equation} \label{eq3.2.3}
\frac{\widetilde{\omega} (z,A,\mathcal{M})}{1-\epsilon}\leq  \widetilde{\omega}(z,A\cap A^{\ast},\mathcal{M}_{\epsilon,A} ),
\qquad z\in\mathcal{M}_{\epsilon,A}.
\end{equation}
To prove the converse inequality of (\ref{eq3.2.3}), let
$u\in\mathcal{PSH}(\mathcal{M}_{\epsilon,A})$ be such that $u\leq 1$ on
$\mathcal{M}_{\epsilon,A}$
and $u\leq 0$ on $ A^{\ast} .$
Consider  the following function
\begin{equation*}
 \hat{u}(z):=
\begin{cases}
\max\left\lbrace (1-\epsilon)u(z),\widetilde{\omega} (z,A,\mathcal{M})\right\rbrace,
  &z\in   \mathcal{M}_{\epsilon,A},\\
  \widetilde{\omega} (z,A,\mathcal{M}), & z \in \mathcal{M}\setminus \mathcal{M}_{\epsilon,A}.
\end{cases}
\end{equation*}
It can be checked that  $\hat{u}\in\mathcal{PSH}(\mathcal{M})$   and $\hat{u}\leq 1.$ Moreover,
in virtue of the assumption on  $u$ and (\ref{eq3.2.2}),   we have that
\begin{equation*} \hat{u}(a)
\leq \max\left\lbrace   (1-\epsilon)u(a),
  \widetilde{\omega}  (a,A,\mathcal{M}) \right\rbrace=0 ,\quad a\in  A^{\ast}  .
\end{equation*}
Consequently,    $\hat{u}\leq\widetilde{\omega} (\cdot,   A^{\ast},\mathcal{M})=\widetilde{\omega} (\cdot,A,\mathcal{M}).$
In particular, one gets that
\begin{equation*}
u(z)\leq \frac{ \widetilde{\omega}(z,A,\mathcal{M})}{1-\epsilon},\qquad z\in \mathcal{M}.
\end{equation*}
Since $u$ is arbitrary, we deduce from the latter estimate that the
converse inequality of (\ref{eq3.2.3}) also holds. This completes the
proof of Part 4).

Part 5) follows immediately from Parts 3) and 4).

Hence, the proof of the proposition is finished.
\end{proof}

The following result shows that our definition of the plurisubharmonic
measure recovers the one given by Alehyane--Zeriahi in \cite[formula
(2.1.2)]{az}.
\begin{prop}\label{prop3.6}
Let $\mathcal{M}$ be a Stein manifold. Let $U$ be a  subdomain of $\mathcal{M}$
 which admits an  exhaustive
sequence of open subsets $\left( U_j \right)_{j=1}^{\infty},$ i.e.
$ U_j \Subset U_{j+1}$ and
$\bigcup\limits_{j=1}^{\infty} U_j= U.$
Then for any  subset $A\subset U,$ there holds
 \begin{equation*}
\widetilde{\omega}(z, A,  U)=\lim\limits_{j\to\infty}
 h^{\ast}_{A\cap U_j, U_j}(z ),\qquad z\in U.
\end{equation*}
\end{prop}
\begin{proof}
First observe that the sequence  $\left(h^{\ast}_{A\cap U_j,
U_j}\right)_{j=1}^{\infty}$ decreases,  as $j\to\infty,$ to  a function $h\in\mathcal{PSH}(\mathcal{M}).$

Next,   since $ U_j\Subset \mathcal{M}$ and $A\setminus A^{\ast}$ is pluripolar, it follows
from Lemma 2.2 in \cite{ah} and Part 1) of Proposition 3.5 that  $h^{\ast}_{A\cap U_j,U_j}=
h^{\ast}_{A\cap A^{\ast}\cap
U_j,U_j}=\widetilde{\omega} (\cdot, A\cap U_j,U_j)$ for any $j\geq 1.$ Consequently, applying Part 2) of Proposition \ref{prop3.5} yields that
 \begin{equation}\label{eq3.2.4}
\widetilde{\omega} (z, A,  U)\leq\liminf\limits_{j\to\infty} \widetilde{\omega}(z, A\cap U_j,U_j)= \lim\limits_{j\to\infty}
 h^{\ast}_{A\cap U_j, U_j}(z )=h(z),\qquad z\in U.
\end{equation}
On the other hand, using the above definition of $h,$ one can check that
$h\leq 1$ on $\mathcal{M}$ and $h\leq 0$ on $\bigcup_{j=1}^{\infty} (A\cap
U_j)^{\ast}.$ Since the latter union is equal to $A^{\ast},$ it follows
that  $h\leq \widetilde{\omega}(\cdot, A,  U).$ This, combined with estimate
(\ref{eq3.2.4}), completes the proof.
\end{proof}

\begin{prop}\label{prop3.7}
Let $\mathcal{M}_j$ be a complex manifold and  $A_j$  a nonempty open
subset of $\mathcal{M}_j,$  $j=1,\ldots,N,$ $N\geq 2.$  \\
1) Then, for $z=(z_1,\ldots,z_N)\in\mathcal{M}_1\times\cdots\times
\mathcal{M}_N,$
 \begin{equation*}
\widetilde{\omega}(z,A_1\times\cdots\times A_N,\mathcal{M}_1\times\cdots\times \mathcal{M}_N)=
\max\limits_{j=1,\ldots,N}\widetilde{\omega}(z_j,A_j,\mathcal{M}_j).
\end{equation*}
2)  Put $\widehat{X}:= \widehat{\X}\left(A_1,\ldots,A_N; \mathcal{M}_1,\ldots, \mathcal{M}_N
\right).$ Then $A_1\times\cdots\times A_N\subset \widehat{X}$ and
\begin{equation*}
\widetilde{\omega}(z,A_1\times\cdots\times A_N, \widehat{X})=
\sum\limits_{j=1}^N\widetilde{\omega}(z_j,A_j,\mathcal{M}_j),\quad
z=(z_1,\ldots,z_N)\in  \widehat{X}.
\end{equation*}
\end{prop}
\begin{proof}
Part 1)  follows immediately by combining Theorem \ref{Rosaythm} and
  the work of  Edigarian and Poletsky in \cite{ep}.

Using Part 1), the proof of Lemma 3 in \cite{jp2} still works in this
context making the obviously necessary changes.
\end{proof}%
\subsection{Three uniqueness theorems and a Two-Constant Theorem}
The following  uniqueness theorems  will play a key role in the sequel.
\begin{thm}\label{unique1}
Let $\mathcal{M}$ be a connected complex manifold,   $A$   a non locally pluripolar subset of $\mathcal{M},$
and    $Z$   a complex analytic space. Let $f,g\in\mathcal{O}(\mathcal{M},Z)$ such that
$f(z)=g(z),$ $z\in A.$  Then $f\equiv g.$
\end{thm}
\begin{proof} Since $A$ is non locally pluripolar, there is an open subset
$U\subset \mathcal{M}$  biholomorphic to an Euclidean domain\footnote{An {\it Euclidean domain} is, by definition,
a domain in $\C^n.$}  such that $A\cap U$ is
nonpluripolar in $U.$ Consequently, we deduce from the equality $f(z)=g(z),$ $z\in A\cap U,$
that $f=g$ on $U.$ Since $\mathcal{M}$ is connected, the desired
conclusion of the theorem follows.
\end{proof}

\begin{thm}  \label{unique2}
  Let  $D_j $ be a complex manifold
  and  $ A_j\subset
D_j$ a non locally pluripolar   subset,
  $j=1,\ldots,N,$  $N\geq 2.$ Let $Y$ be a complex analytic  space. 
Let $U_1$ and $U_2$ be  two open subsets of $D_1.$
For $k\in\{1,2\}, $  let $f_k\in\mathcal{O}(\widehat{X}_k,Y),$  where
\begin{equation*}
\widehat{X}_k:=\widehat{\X}\left(A_1\cap U_k,A_2,\ldots,A_N;U_k,D_2,\ldots,D_N\right).
\end{equation*}
Then:\\
1) If
$ f_1=f_2$ on $(U_1\cap  U_2)\times (A_2\cap A^{\ast}_2)\times\cdots\times (A_N\cap A^{\ast}_N) ,$
 then $f_1=f_2$ on $\widehat{X}_1\cap\widehat{X}_2.$\\
 2) If $U_1=U_2$ and $ f_1=f_2$ on $(A_1\cap A^{\ast}_1\cap U_1 )\times (A_2\cap A^{\ast}_2)
 \times\cdots\times (A_N\cap A^{\ast}_N) ,$
 then  $f_1=f_2$ on $\widehat{X}_1 .$
\end{thm}
 \begin{proof} To prove Part 1), fix an arbitrary point  $z^0=(z_1^0,\ldots,z_N^0)\in\widehat{X}_1
\cap\widehat{X}_2.$  We need to show that $ f_1(z^0)=f_2(z^0).$

For any $2\leq j\leq N,$ let $\mathcal{G}_j$ be the connected component which contains $z_j^0$ of the
following open set
\begin{equation*}
  \left\lbrace z_j\in D_j:\ \widetilde{\omega} (z_j,A_j,D_j)<1-\max\limits_{k\in\{1,2\}} \widetilde{\omega}(z_1^0,A\cap
  U_{k}, U_k)
  -\sum\limits_{p=2}^{j-1} \widetilde{\omega}(z^0_p,A_p,D_p) \right\rbrace.
  \end{equation*}

 Observe that for $k\in\{1,2\}$ and $(a_3,\ldots,a_N)\in
  (A_3\cap A^{\ast}_3)\times\cdots\times (A_N\cap A^{\ast}_N) ,$
 the mapping  $z_2\in\mathcal{G}_2\mapsto f_k\left(z^0_1,z_2,a_3,\ldots,a_N\right)$
 belongs to  $\mathcal{O}
(\mathcal{G}_2,Y).$ In addition, it follows from  the hypothesis that
\begin{equation}\label{eq3.3.1}
  f_1\left(z^0_1,a_2,\ldots,a_N\right)= f_2\left(z^0_1,a_2,\ldots,a_N\right),
  \qquad  a_2 \in A_2\cap A^{\ast}_2  .
  \end{equation}
On the other hand, by Part 5) of Proposition \ref{prop3.5},   $\mathcal{G}_2$ contains a
non locally pluripolar subset of $A_2\cap A^{\ast}_2.$ Therefore, by Theorem  \ref{unique1},
\begin{multline*}
 f_1\left(z^0_1,z_2,a_3,\ldots,a_N\right)= f_2\left(z^0_1,z_2,a_3,\ldots,a_N\right),\\
  \qquad \left(z_2,a_3,\ldots,a_N\right)
  \in \mathcal{G}_2\times (A_3\cap A^{\ast}_3)\times\cdots\times(A_N\cap A^{\ast}_N) .
\end{multline*}
Hence,
\begin{multline}\label{eq3.3.2}
 f_1\left(z^0_1,z_2^0,a_3,\ldots,a_N\right)=
 f_2\left(z^0_1,z_2^0,a_3,\ldots,a_N\right),\\
  \qquad \left(a_3,\ldots,a_N\right)
  \in  (A_3\cap A^{\ast}_3)\times\cdots\times(A_N\cap A^{\ast}_N) .
\end{multline}

Repeating the  argument in (\ref{eq3.3.1})--(\ref{eq3.3.2}) $(N-2)$ times, we finally obtain
 $f_1(z^0)=f_2(z^0).$ Hence, the proof of  Part 1) is finished.

In  virtue of Part 1),  Part 2) is reduced to proving that
\begin{equation*}
f_1=f_2\qquad\text{on}\  U_1\times (A_2\cap A^{\ast}_2)\times\cdots\times( A_N\cap A^{\ast}_N).
\end{equation*}
To do this, fix an
arbitrary point
  $z^0=(z_1^0,a^0_2\ldots,a_N^0)\in U_1\times (A_2\cap A^{\ast}_2)\times\cdots\times( A_N\cap A^{\ast}_N)  $
  such that $z^0\in\widehat{X}_1.$
 Then $ \widetilde{\omega} (z_1,A_1\cap A^{\ast}_1\cap U_1,U_1)< 1 .$ Let $\mathcal{G}$  be the connected  component
 containing $z^0_1$ of $U_1.$ Using Parts 1) and 3)  of Proposition
 \ref{prop3.5} and taking into account the latter estimate, we see that
$A_1\cap A^{\ast}_1\cap\mathcal{G}$ is a non locally pluripolar set.

 Next, observe that for $k\in\{1,2\},$
 the mapping  $z_1\in\mathcal{G}\mapsto f_k\left(z_1,a^0_2,\ldots,a^0_N\right)$
 belongs to  $\mathcal{O}
(U_1,Y).$
 Moreover, since we know from the hypothesis and the above paragraph that  $ f_1\left(\cdot,a^0_2,\ldots,a^0_N\right)=f_2\left(\cdot,a^0_2,\ldots,a^0_N\right)$
  on the non locally pluripolar set $ A_1\cap A^{\ast}_1\cap \mathcal{G}
,$ it follows from  Theorem  \ref{unique1}  that
\begin{equation*}
  f_1\left(z_1,a^0_2,\ldots,a^0_N\right)= f_2\left(z_1,a^0_2,\ldots,a^0_N\right),
  \qquad  z_1
  \in   \mathcal{G} .
  \end{equation*}
Hence, $f_1(z^0)=f_2(z^0),$ which completes the proof of Part 2).
\end{proof}

The next result, combined with Part 2) of Theorem \ref{unique2}, establishes the uniqueness stated in Theorem A and in    its intermediate versions (see Theorems
\ref{thm4.1}, \ref{thm5.1}, \ref{thmA_special_case} and Proposition \ref{prop_special_case}  below).

 \begin{thm}\label{unique3}
Let  $D_j $ be a complex manifold
  and  $ A_j\subset
D_j$ a non locally pluripolar   subset,
  $j=1,\ldots,N.$ Let $Z$ be a complex analytic space. One defines $X,$ $X^{\ast}$ and $\widehat{X}$
  as in Subsection 2.1. Let
  $f\in\mathcal{O}_s(X,Z)$ and
$\hat{f}\in\mathcal{O}(\widehat{X},Z)$
such that $\hat{f}=f$ on $X\cap X^{\ast}.$
Then  $\hat{f}=f$ on
$X\cap \widehat{X}.$
\end{thm}
\begin{proof}
Let $z^0=(z^0_1,\ldots,z^0_N)$ be an arbitrary point of $X\cap
\widehat{X},$ and put $f_1:=\hat{f},$  $f_2:=f.$ Arguing as in the proof
of Theorem \ref{unique2}, we can show that $\hat{f}(z^0)=f(z^0).$
This completes the proof.
\end{proof}

\smallskip
The following Two-Constant Theorem for plurisubharmonic  functions
will play an important  role in  the proof of the estimate in Theorem A.
\begin{thm}\label{two-constant}
Let $\mathcal{M}$ be a   complex manifold  and $A$   a non locally pluripolar subset of $\mathcal{M}.$
   Let  $m,M\in\R$  and   $u\in\mathcal{PSH}(\mathcal{M})$ such that
$u(z)\leq M$ for $z\in\mathcal{M},$  and $u(z)\leq m$ for $z\in
A.$ 
Then
\begin{equation*}
u(z)\leq m(1- \widetilde{\omega}(z,A,\mathcal{M}))+M\cdot \widetilde{\omega}(z,A,\mathcal{M}),\qquad z\in  \mathcal{M}.
\end{equation*}
\end{thm}
\begin{proof}
It follows immediately from the definition of $\widetilde{\omega}(\cdot,A,\mathcal{M})$ given in
Subsection 2.1.
\end{proof}
\section{Part 1 of the proof of Theorem A }

The main purpose of the section is  to prove Theorem A in the following
special case.
\begin{thm}\label{thm4.1}
Let  $D$ be a  complex manifold, let  $G$ be  a complex manifold which is biholomorphic to an open set in $\C^q$
 $(q\in\N),$
       let $ A$    be an open  subset of $D,$  and let $B $ be a non locally pluripolar subset of $G.$ Let $Z$
        be a complex  analytic space possessing the Hartogs
  extension property. Put $X:=\X(A,B;D,G)$ and
$\widehat{X}:=\widehat{\X}(A,B;D,G).$
Then for any mapping $f\in\mathcal{O}_s(X,Z),$ there is a unique mapping
$\hat{f}\in\mathcal{O}(\widehat{X},Z)$
such that $\hat{f}=f$ on $X\cap X^{\ast} .$
\end{thm}
\begin{rem}\label{rem4.2}
Under the above hypothesis, it can be checked that $X\cap X^{\ast}=(A\times G)\cup (D\times(B\cap B^{\ast})).$ In
addition, in the proof below  we  assume that   $G$ is a domain in $\C^q.$
Clearly, in virtue of Part 3) of Proposition \ref{prop3.5}, it suffices to prove the theorem  under this assumption.
\end{rem}
\begin{proof}
 We begin the proof with the following lemma.
\begin{lem}
\label{lem4.3}
We keep the hypothesis of Theorem \ref{thm4.1}. For $j\in \{1,2\},$ let $\phi_j\in\mathcal{O}(
 \overline{E}, D) $ be a holomorphic disc,  and let $t_j\in E$ such that $\phi_1(t_1)=\phi_2(t_2)$ and
$  \frac{1}{2\pi}\int\limits_{0}^{2\pi}1_{D\setminus A} (\phi_j(e^{i\theta}))d\theta
<1.$ Then:
\\
1) For  $j \in \{
1,2\},$  the mapping $(t,w)\mapsto f(\phi(t),w)$ belongs to
$ \mathcal{O}_s(\X(\phi^{-1}_j(A)\cap E,B;E,G),Z), $ where
$\phi^{-1}_j(A):=\lbrace t\in \overline{E}:\ \phi_j(t)\in A\rbrace.$     \\
2)  For  $j \in \{
1,2\},$ in virtue of Part 1), Remark \ref{rem4.2} and applying Theorem 3,
 let $\hat{f}_j$ be the unique mapping in
$ \mathcal{O}\left(\widehat{\X}(\phi^{-1}_j(A)\cap E,B;E,G),Z\right)$ such that  $\hat{f}_j(t,w)=f(\phi_j(t),w),$
 $(t,w)\in\X\left(\phi^{-1}_j(A)\cap E,B\cap B^{\ast};E,G\right).$   Then
\begin{equation*}
\hat{f}_1(t_1,w)=\hat{f}_2(t_2,w),
\end{equation*}
for all $w\in G$ such that
$(t_j,w)\in\widehat{\X}\left(\phi^{-1}_j(A)\cap E,B;E,G\right),$  $j \in \{1,2\}.$
\end{lem}

\smallskip
\noindent {\it Proof of Lemma \ref{lem4.3}.}
 Part 1) follows immediately from the hypothesis. Therefore, it remains to prove Part
 2). To do this fix  $w_0\in G$ such that $(t_j,w_0)\in\widehat{\X}\left(\phi^{-1}_j(A)\cap E
 ,B;E,G\right)$
 for  $j \in \{1,2\}.$
We need to show that $\hat{f}_1(t_1,w_0)=\hat{f}_2(t_2,w_0).$  Observe
that both mappings  $w\in\mathcal{G}\mapsto \hat{f}_1(t_1,w)$ and $w\in\mathcal{G}\mapsto\hat{f}_2(t_2,w)$ belong to  $\mathcal{O}
(\mathcal{G},Z),$
where $\mathcal{G}$ is the connected component which contains $w_0$ of the
following open set
\begin{equation*}
  \left\lbrace w\in G:\ \widetilde{\omega}(w,B,G)<1-\max\limits_{ j \in \{1,2\}}  \widetilde{\omega}(t_j,\phi^{-1}_j(A)\cap E,E)
   \right\rbrace.
   \end{equation*}
 Moreover, for any $w\in B\cap B^{\ast}$ we have that   $\widetilde{\omega}(w, B\cap B^{\ast},G)=0.$
 Therefore,    for any  $j \in \{
1,2\}$  one clearly gets that
   $(t_j,w)\in\X\left(\phi^{-1}_j(A)\cap E,B\cap B^{\ast};E,G\right).$   Consequently, since
   $\phi_1(t_1)=\phi_2(t_2),$ it follows that
\begin{equation*}
\hat{f}_1(t_1,w)=f(\phi_1(t_1),w)=f(\phi_2(t_2),w)=\hat{f}_2(t_2,w) , \qquad w\in B\cap B^{\ast}.
\end{equation*}
On the other hand, by Part 5) of Proposition \ref{prop3.5},   $\mathcal{G}$ contains a
non locally pluripolar subset of $B\cap B^{\ast}.$ Therefore, by Theorem  \ref{unique1},   $\hat{f}_1(t_1,w)=\hat{f}_2(t_2,w),$
  $  w\in\mathcal{G}.$ Hence,  $\hat{f}_1(t_1,w_0)=\hat{f}_2(t_2,w_0),$
which completes the proof of the lemma.
\hfill $\square$

\smallskip

\noindent{\bf Step 1:} {\it Construction of the extension  mapping $\hat{f}$ on $\widehat{X}.$}

\smallskip

\noindent{ \it Proof of Step 1.} We define $\hat{f}$ as follows: Let $\mathcal{X}$ be  the set of all pairs
$(z,w)\in D\times G$  with the property that there are a holomorphic disc
$\phi\in\mathcal{O}(\overline{E},D)$ and $t\in E$ such that $\phi(t)=z$
and  $(t,w)\in\widehat{\X}\left(\phi^{-1}(A)\cap E,B;E,G\right).$  In virtue of
 Theorem 3,
 let $\hat{f}_{\phi}$ be the unique mapping  in
$ \mathcal{O}\left(\widehat{\X}(\phi^{-1}(A)\cap E,B;E,G),Z\right)$ such that
\begin{equation}\label{eq4.1}
  \hat{f}_{\phi}(t,w)=f(\phi(t),w),\qquad (t,w)\in \X\left(\phi^{-1}(A)\cap E,B\cap B^{\ast};E,G\right).
\end{equation}
   Then  the desired extension mapping $\hat{f}$ is given by
\begin{equation}\label{eq4.2}
 \hat{f}(z,w):=\hat{f}_{\phi}(t,w) .
\end{equation}
In virtue of Part 2) of Lemma \ref{lem4.3}, $\hat{f}$ is well-defined on $\mathcal{X}.$
We next prove that
\begin{equation}
 \label{eq4.3}
 \mathcal{X}=\widehat{X}  .
\end{equation}
Taking  (\ref{eq4.3}) for granted, then $\hat{f}$ is
well-defined on $\widehat{X}.$ Moreover, it follows from formula
(\ref{eq4.2}) that for every
fixed $z\in D,$  the restricted mapping $\hat{f}(z,\cdot)$ is holomorphic on the open set
$\left\lbrace w\in G:\ (z,w)\in \widehat{X} \right\rbrace.$

Now we return to (\ref{eq4.3}). To prove the inclusion $\mathcal{X}\subset\widehat{X},$
let  $(z,w)\in\mathcal{X}.$  By the above definition of $\mathcal{X},$ one may find a  holomorphic  disc    $\phi\in
\mathcal{O}(\overline{E}, D),$  a point  $t\in E$ such that $\phi(t)=z$
and  $(t,w)\in\widehat{\X}\left(\phi^{-1}(A)\cap E,B;E,G\right).$
Since  $ \widetilde{\omega}(\phi(t),A,D)\leq  \widetilde{\omega}(t,\phi^{-1}(A)\cap E,E),$ it follows that
\begin{equation*}
 \widetilde{\omega} (z,A,D) +  \widetilde{\omega}(w,B,G) \leq   \widetilde{\omega}(t,\phi^{-1}(A)\cap E,E)+
  \widetilde{\omega}(w,B,G)<1,
\end{equation*}
Hence  $(z,w)\in\widehat{X}.$ This proves the above mentioned inclusion.

To finish the proof of  (\ref{eq4.3}), it suffices to show that  $\widehat{X}\subset\mathcal{X}.$
To do this, let $(z,w)\in\widehat{X}$ and fix any $\epsilon>0$ such that
\begin{equation}\label{eq4.4}
\epsilon<1-  \widetilde{\omega}(z,A,D) -  \widetilde{\omega}(w,B,G).
\end{equation}
 Applying  Theorem \ref{Rosaythm} and Proposition \ref{prop3.4}, there is a holomorphic  disc   $\phi\in
\mathcal{O}(\overline{E},  D)$   such that $\phi(0)=z$ and
\begin{equation}\label{eq4.5}
\frac{1}{2\pi} \int\limits_{0}^{2\pi} 1_{D\setminus
A}(\phi(e^{i\theta}))d\theta< \widetilde{\omega}(z,A,D)+\epsilon.
\end{equation}
Observe that
\begin{eqnarray*}
  \widetilde{\omega}(0,\phi^{-1}(A)\cap E,E)+ \widetilde{\omega}(w,B,G)&\leq &\frac{1}{2\pi} \int\limits_{0}^{2\pi} 1_{D\setminus
A}(\phi(e^{i\theta}))d\theta+ \widetilde{\omega}(w,B,G)\\
 &< &  \widetilde{\omega}(z,A,D) +  \widetilde{\omega}(w,B,G)+\epsilon<1,
\end{eqnarray*}
where the first inequality follows from an application of  Lemma \ref{lem3.3},
the second one from (\ref{eq4.5}), and the last one from (\ref{eq4.4}). Hence,
 $(0,w)\in\widehat{\X}\left(\phi^{-1}(A)\cap E,B;E,G\right),$ which
 implies that $(z,w)\in\mathcal{X}.$ This complete  the proof of  (\ref{eq4.3}).
Hence Step I is finished. \hfill $\square$

Next, we would like   to show that $\hat{f}$ satisfies the conclusion of the
theorem.  This will be accomplished  in two steps below.

\smallskip

\noindent{\bf Step 2:} {\it  Proof of the equality $\hat{f}=f$ on $X\cap X^{\ast}.$}

\smallskip

\noindent{\it Proof of Step 2.}
Let $(z,w)$ be an arbitrary point of $ A\times G.$
   Choose the   holomorphic  disc
$\phi\in \mathcal{O}(\overline{E},D)$ given by $\phi(t):=z,$ $
t\in\overline{E}.$ Then  by formula (\ref{eq4.2}),
\begin{equation}\label{eq4.6}
\hat{f}(z,w)=\hat{f}_{\phi}(0,w)=f(\phi(0),w)=f(z,w),  \qquad w\in G.
\end{equation}
Hence, $\hat{f}=f$ on $A\times G.$

Next, let $(z,w)$ be an arbitrary point of $ D\times (B\cap B^{\ast})$ and  let
$\epsilon>0$ be
such that
\begin{equation}\label{eq4.7}
\epsilon<1- \widetilde{\omega}(z,A,D).
\end{equation}
 Applying Theorem \ref{Rosaythm} and Proposition
\ref{prop3.4}, one may find  a holomorphic disc
$\phi\in\mathcal{O}( \overline{E} , D)$  such that  $\phi(0)=z$ and
\begin{equation}\label{eq4.8}
\frac{1}{2\pi} \int\limits_{0}^{2\pi} 1_{D\setminus
A}(\phi(e^{i\theta}))d\theta< \widetilde{\omega}(z,A,D)+\epsilon.
\end{equation}
Consequently,
\begin{equation*}
 \widetilde{\omega}(0,\phi^{-1}(A)\cap E,E)+ \widetilde{\omega}(w,B,G)\leq \frac{1}{2\pi} \int\limits_{0}^{2\pi} 1_{D\setminus
A}(\phi(e^{i\theta}))d\theta< \widetilde{\omega}(z,A,D)+\epsilon<1,
\end{equation*}
where the first inequality follows from an application of  Lemma \ref{lem3.3}
and the equality $ \widetilde{\omega}(w,B,G)=0,$
the second one from (\ref{eq4.8}), and the last one from (\ref{eq4.7}).
Hence,  $(0,w)\in\widehat{\X}\left(\phi^{-1}(A)\cap E,B;E,G\right).$
Therefore,  using  (\ref{eq4.1})--(\ref{eq4.2}) and arguing as in (\ref{eq4.6}),
we conclude  that $\hat{f}(z,w) = f(z,w).$ This proves that $\hat{f}=f$ on
$D\times ( B\cap B^{\ast}).$

In summary, we have shown that  $\hat{f}=f$ on
$(A\times G)\cup \left(D\times ( B\cap B^{\ast})\right).$ In virtue of Remark  \ref{rem4.2}, Step 2 is complete. \hfill $\square$

\smallskip

\noindent{\bf Step 3:} {\it  Proof of the fact that
$\hat{f}\in\mathcal{O}(\widehat{X},Z).$}

\smallskip

\noindent{\it Proof of Step 3.} Fix an arbitrary point   $(z_0,w_0)\in
\widehat{X}$ and let  $\epsilon>0$ be so small  such that
\begin{equation}\label{eq4.9}
2\epsilon<1-  \widetilde{\omega}(z_0,A,D) -  \widetilde{\omega}(w_0,B,G).
\end{equation}
Since $ \widetilde{\omega}(\cdot,B,G)\in\mathcal{PSH}(G),$   one may find an open neighborhood $V $ of $w_0
$ such that
\begin{equation}\label{eq4.10}
   \widetilde{\omega}(w,B,D) <  \widetilde{\omega}(w_0,B,G)+\epsilon,\qquad  w\in V.
\end{equation}
Let $d$ be the dimension of $D$ at the point $z_0.$ Applying Lemma
\ref{Rosaylem} and Proposition \ref{prop3.4}, we obtain an open set $T$ in $\C,$ an open  neighborhood $U$ of $z_0$
which is biholomorphic to the unit ball in  $\C^d,$ and
  and a  family of holomorphic discs $(\phi_z)_{z\in U}\subset \mathcal{O}(\overline{E},D)$ with the
following properties:
\begin{eqnarray} \label{eq4.11}
\text{the mapping}\ (z,t)\in U\times E &\mapsto& \phi_z(t)\ \text{is
holomorphic};  \\\label{eq4.12}
\phi_{z}(0)&=&z,\qquad z\in U;\\
\label{eq4.13}\phi_z(t)&\in & A, \qquad t\in T\cap\overline{E},\ z\in U;\\
  \label{eq4.14}\frac{1}{2\pi}\int\limits_{0}^{2\pi}1_{\partial E\setminus T  } (e^{i\theta})d\theta
&<&  \widetilde{\omega}(z_0,A,D)+\epsilon.
\end{eqnarray}

Consider the   mapping $g:\ \X\left(T\cap E,U,B;E,U, G\right)\longrightarrow Z$ given  by
\begin{equation}\label{eq4.15}
g(t,z,w):=f(\phi_z(t),w),\qquad (t,z,w)\in\X\left(T\cap E,U,B;E,U, G\right).
\end{equation}
We make the following observations:

  Let $t\in T\cap E.$   Then, in virtue of (\ref{eq4.13}) we have $\phi_z(t)\in
 A$ for  $z\in U.$ Consequently, in virtue of  (\ref{eq4.11}), (\ref{eq4.15}) and the hypothesis
 $f\in\mathcal{O}_s(X,Z),$ we conclude that
 $g(t,z,\cdot)|_G\in\mathcal{O}(G,Z)$ $\Big($resp.
 $g(t,\cdot,w)|_U\in\mathcal{O}(U,Z)\Big)$ for any $z\in U$ (resp. $w\in B).$
Analogously, for any $z\in U,$  $w\in B,$    we can show that
 $g(\cdot,z,w)|_E\in\mathcal{O}(E,Z).$

 In summary, we have shown that  $g\in \mathcal{O}_s\left(
 \X\left(T\cap E,U,B;E,U, G\right), Z\right).$ Recall that  $U$ is biholomorphic to the unit ball in $\C^d.$
  Consequently, we are able to apply
 Theorem 3 to $g$ in order to obtain a unique mapping $\hat{g}\in
\mathcal{O}\left(
 \widehat{\X}\left(T\cap E,U,B;E,U, G\right), Z\right)$ such that
\begin{equation}\label{eq4.16}
\hat{g}(t,z,w)= g(t,z,w),\qquad (t,z,w)\in\X\left(T\cap E,U,B\cap B^{\ast};E,U, G\right).
\end{equation}
Observe that
\begin{equation*}
  \widehat{\X}\left(T\cap E,U,B;E,U, G\right)=\left\lbrace  (t,z,w)\in E\times U\times G:\
   \widetilde{\omega}(t,T\cap E,E)+ \widetilde{\omega}(w,B,G)<1\right\rbrace.
\end{equation*}
On the other hand, for any $w\in V,$
\begin{equation}\label{eq4.17}
\begin{split}
  \widetilde{\omega}(0,T\cap E,E)+ \widetilde{\omega}(w,B,G)&\leq \frac{1}{2\pi} \int\limits_{0}^{2\pi}
 1_{\partial E\setminus T}(e^{i\theta})d\theta+ \widetilde{\omega}(w_0,B,G)+\epsilon\\
 &<   \widetilde{\omega}(z_0,A,D) +  \widetilde{\omega}(w_0,B,G)+2\epsilon<1,
 \end{split}
\end{equation}
where the first inequality follows from an application of  Lemma \ref{lem3.3} and
(\ref{eq4.10}),
the second one from (\ref{eq4.14}), and the last one from (\ref{eq4.9}).
 Consequently,
\begin{equation}\label{eq4.18}
 (0,z,w)\in\widehat{\X}\left(T\cap E,U,B;E,U, G\right),\qquad (z,w)\in U\times V.
\end{equation}
 It  follows from   (\ref{eq4.2}),  (\ref{eq4.12}),  (\ref{eq4.13}) and  (\ref{eq4.17}) that, for $z \in U ,$
  $\hat{f}_{\phi_z}$ is well-defined and holomorphic on $\widehat{\X}(T\cap E,B;E,
G),$ and
\begin{equation}\label{eq4.19}
 \hat{f}(z,w)=\hat{f}_{\phi_z}(0,w),\qquad w\in V.
\end{equation}
On the other hand, it follows from  (\ref{eq4.1}), (\ref{eq4.15}) and  (\ref{eq4.16}) that
\begin{equation*}
\hat{f}_{\phi_z}(t,w)=\hat{g}(t,z,w),\qquad (t,w)\in\X\left(T\cap E,B\cap B^{\ast};E, G\right), z \in U .
\end{equation*}
 Since, for fixed $z\in U,$ the restricted mapping
 $(t,w)\mapsto \hat{g}(t,z,w)$ is holomorphic on $\widehat{\X}(T\cap E,B;E, G),$
 we deduce from the latter equality and the uniqueness of Theorem 3 that
\begin{equation*}
\hat{g}(t,z,w)=\hat{f}_{\phi_z}(t,w),\qquad (t,w)\in\widehat{\X}\left(T\cap E,B\cap B^{\ast};E, G\right), z \in U .
\end{equation*}
In particular, using  (\ref{eq4.2}), (\ref{eq4.18}) and (\ref{eq4.19}),
\begin{equation*}
\hat{g}(0,z,w)=\hat{f}_{\phi_z}(0,w)=\hat{f}(z,w),\qquad   (z,w) \in U\times V .
\end{equation*}
Since we know from (\ref{eq4.18}) that $\hat{g}$ is holomorphic on a neighborhood of $(0,z_0,w_0),$
we conclude that $\hat{f}$ is holomorphic on a neighborhood of $(z_0,w_0).$
 Since $(z_0,w_0)\in\widehat{X}$  is arbitrary,
it follows that  $\hat{f}\in\mathcal{O}(\widehat{X},Z).$ Hence Step 3 is complete.
 \hfill $\square$

 Combining Steps 1--3, the theorem follows.
\end{proof}

\section{Part 2 of the proof of Theorem A }

The main purpose of the section is  to prove Theorem A in the following
special case.
\begin{thm}\label{thm5.1}
Let  $D,$ $G$ be  complex manifolds,    and let $ A\subset
D,$  $B\subset G$  be open  subsets.  Let $Z$ be a complex analytic space possessing the Hartogs extension
property.  Put $X:=\X(A,B;D,G)$ and
$\widehat{X}:=\widehat{\X}(A,B;D,G).$
Then for any mapping $f\in\mathcal{O}_s(X,Z),$ there is a unique   mapping
$\hat{f}\in\mathcal{O}(\widehat{X},Z)$
such that $\hat{f}=f$ on $X.$
\end{thm}
\begin{rem}\label{rem5.2} Using Part 1) of Proposition \ref{prop3.5}
 it  can be checked that, under the above hypothesis,
$X^{\ast}=X.$
\end{rem}
\begin{proof}  Let us start  with the following lemma.
\begin{lem}
\label{lem5.3}
We keep the hypothesis of Theorem 5.1. For $j\in\{ 1,2\},$ let $\psi_j \in\mathcal{O}( \overline{E},
  G)$
be a holomorphic disc and let $\tau_j\in E$ such that $\psi_1(\tau_1)=\psi_2(\tau_2)$ and
$  \frac{1}{2\pi}\int\limits_{0}^{2\pi} 1_{G\setminus B}(\psi_j(e^{i\theta})) d\theta
<1.$ Then:
\\
1) For $j\in\{ 1,2\},$ the mapping
$(z,\tau) \mapsto f(z,\psi_j(\tau))$ belongs to $ \mathcal{O}_s\left(\X\left(A,\psi^{-1}_j(B)\cap E;D,E
\right),Z\right), $ where
$\psi^{-1}_j(B):=\lbrace \tau\in \overline{E}:\ \psi_j(\tau)\in B\rbrace.$     \\
2)  For $j\in\{ 1,2\},$ in virtue of Part 1), Remark \ref{rem5.2} and applying  Theorem \ref{thm4.1},
   let $\tilde{f}_j$ be the unique  mapping in
$ \mathcal{O}\left(\widehat{\X}\left(A,\psi^{-1}_j(B)\cap E;D,E\right),Z\right)$ such that
$\tilde{f}_j(z,\tau)=f(z,\psi_j(\tau)),$
 $(z,\tau)\in\X\left(A,\psi^{-1}_j(B)\cap E;D,E\right).$  Then
\begin{equation*}
\tilde{f}_1(z,\tau_1)=\tilde{f}_2(z,\tau_2),
\end{equation*}
for all $z\in D$ such that
$(z,\tau_j)\in\widehat{\X}\left(A,\psi^{-1}_j(B)\cap E;D,E\right),$   for $j\in\{
1,2\}.$
\end{lem}

\smallskip
\noindent {\it Proof of Lemma \ref{lem5.3}.}
It follows along the same lines as those of Lemma \ref{lem4.3}.
\hfill $\square$

\smallskip

Now we return to Theorem  \ref{thm5.1}.
First we define  a mapping $\hat{g}:\ \widehat{X}\longrightarrow Z$ as follows: Let $\mathcal{X}$ be  the set of all pairs
$(z,w)\in D\times G$  with the property that there are a holomorphic disc
$\psi\in\mathcal{O}(\overline{E}, G)$ and $\tau\in E$ such that $\psi(\tau)=w$
and  $(z,\tau)\in\widehat{\X}\left(A,\psi^{-1}(B)\cap E;D,E\right).$  In virtue of
Lemma \ref{lem5.3},
 let $\tilde{f}_{\psi}$ be the unique  mapping in
$ \mathcal{O}\left(\widehat{\X}\left(A,\psi^{-1}(B)\cap E;D,E\right),Z\right)$ such that
\begin{equation}\label{eq5.1}
\tilde{f}_{\psi}(z,\tau)=f(z,\psi(\tau)),\qquad  (z,\tau)\in
 \X\left(A,\psi^{-1}(B)\cap E;D,E\right).
 \end{equation}
  Then we define
\begin{equation}
\label{eq5.2} \hat{g}(z,w):=\tilde{f}_{\psi}(z,\tau) .
\end{equation}
In virtue of Part 2) of  Lemma \ref{lem5.3}, $\hat{g}$ is well-defined on $\mathcal{X}.$
Moreover, arguing as in Step 1 in the proof of Theorem 4.1, we  can show
that  $\mathcal{X}=\widehat{X}  .$
Consequently, $\hat{g}$ is
well-defined on $\widehat{X}.$ Moreover, arguing as in Step 2 in the proof
of Theorem \ref{thm4.1}, it follows from (\ref{eq5.1})--(\ref{eq5.2}) and Remark \ref{rem5.2} that for every fixed $w\in G,$ the restricted mapping
 $\hat{g}(\cdot,w)$ is holomorphic  on the open set  $\left\lbrace z\in D:\ (z,w)\in\widehat{X}
 \right\rbrace$  and that
$\hat{g}=f$ on $X.$

We define the desired mapping $\hat{f}$ on $\widehat{X}$ as follows: Let
$(z,w)\in \widehat{X},$   let  $\phi\in
\mathcal{O}(\overline{E} ,D)$ be a  holomorphic disc and $t\in E$ such that $\phi(t)=z$
and  $( t,w)\in\widehat{\X}\left(\phi^{-1}(A)\cap E,B;E,G\right).$  In virtue of
Lemma \ref{lem5.3} and replacing the role of $B$ (resp. $D$) by that of $A$ (resp. $G$) therein,
 let $\hat{f}_{\phi}$ be the unique  mapping in
$ \mathcal{O}\left(\widehat{\X}\left(\phi^{-1}(A)\cap E,B;E,G\right),Z\right)$ such that
\begin{equation}\label{eq5.3}
\hat{f}_{\phi}(t,w)=f(\phi(t), w),\qquad
 (t,w)\in \X(\phi^{-1}(A)\cap E,B;E,G).
\end{equation}
 Then we put
\begin{equation}
\label{eq5.4} \hat{f}(z,w):=\hat{f}_{\phi}(t,w) .
\end{equation}
Arguing as in the previous paragraph, we conclude that $\hat{f}$ is well-defined on $\widehat{X}.$
  Moreover, it follows from (\ref{eq5.3})--(\ref{eq5.4}) and Remark
  \ref{rem5.2}
that  for every fixed $z\in D,$ the restricted mapping
 $\hat{f}(z,\cdot)$ is holomorphic  on the open set  $\left\lbrace w\in G:\ (z,w)\in\widehat{X}
 \right\rbrace$  and that
$\hat{f}=f$ on $X.$

The proof of the theorem will be complete if we can show that
\begin{equation}\label{eq5.5}
\hat{f}=\hat{g}.
\end{equation}
 Indeed, taking (\ref{eq5.5}) for granted, then for any $(z_0,w_0)\in \widehat{X},$
 we may find an open
  neighborhood $U\times V$ of $(z_0,w_0)$ such that $U\times V\subset
  \widehat{X}$ and $U$ (resp. $V$) is biholomorphic to an Euclidean ball.
Using  (\ref{eq5.5}) and  the above-mentioned property of $\hat{f}$ and $\hat{g},$
we see that  $\hat{f} (=\hat{g})\in \mathcal{O}_s(\X(U,V;U,V),Z).$ Consequently, applying Theorem 3
  to $\hat{f} ,$  it
  follows that $\hat{f}\in\mathcal{O}(U\times V,Z).$ Hence,
  $\hat{f}\in\mathcal{O}(\widehat{X},Z),$ and the proof of the theorem is   finished.

To prove (\ref{eq5.5}), fix an arbitrary point $(z_0,w_0)\in\widehat{X}.$ Fix any $\epsilon>0$ such that
\begin{equation}\label{eq5.6}
3\epsilon<1- \widetilde{\omega}(z_0,A,D) -  \widetilde{\omega}(w_0,B,G).
\end{equation}
 Applying Theorem \ref{Rosaythm}  and Proposition \ref{prop3.4}, there is a holomorphic disc   $\phi\in
\mathcal{O}(\overline{E} , D)$  (resp.  $\psi\in\mathcal{O}(
\overline{E} , G))$  such that $\phi(0)=z_0$  (resp. $\psi(0)=w_0$) and
\begin{equation*}
\frac{1}{2\pi} \int\limits_{0}^{2\pi} 1_{D\setminus
A}(\phi(e^{i\theta}))d\theta< \widetilde{\omega}(z_0,A,D)+\epsilon,\qquad
\frac{1}{2\pi} \int\limits_{0}^{2\pi} 1_{G\setminus
B}(\psi(e^{i\theta}))d\theta< \widetilde{\omega}(w_0,B,G)+\epsilon.
\end{equation*}
Using this and estimate (\ref{eq5.6}), and arguing as in Step 1 of Theorem \ref{thm4.1}, we see that
  $(0,0)\in\widehat{\X}\left(\phi^{-1}(A)\cap E,\psi^{-1}(B)\cap E;E,E\right)
 .$ Moreover, since $f\in\mathcal{O}_s(X,Z),$ the mapping $h$ given by
 \begin{equation*}
  h(t,\tau):= f(\phi(t),\psi(\tau)),\qquad (t,\tau)\in
  \X\left(\phi^{-1}(A)\cap E,\psi^{-1}(B)\cap E;E,E\right),
 \end{equation*}
 belongs to $\mathcal{O}_s\Big(\X\left(\phi^{-1}(A)\cap E,\psi^{-1}(B)\cap E;E,E\right),Z\Big).$ Moreover,
 in virtue of (\ref{eq5.1}) and (\ref{eq5.3}),
\begin{equation}\label{eq5.7}
\hat{f}_{\phi}(t,\psi(\tau))=f(\phi(t),\psi(\tau))=\tilde{f}_{\psi}(\phi(t),\tau),\
\ (t,\tau)\in\X\left(\phi^{-1}(A)\cap E,\psi^{-1}(B)\cap E;E,E\right).
\end{equation}
  By Theorem  3, let $\hat{h}\in \mathcal{O}\left(\widehat{\X}\left(\phi^{-1}(A)\cap E,\psi^{-1}(B)\cap E;E,E
 \right),Z\right)$  be the unique  mapping  such that
 \begin{equation*}
 \hat{h}(t,\tau)=h(t,\tau)= f(\phi(t),\psi(\tau)),\qquad
(t,\tau)\in \X\left(\phi^{-1}(A)\cap E,\psi^{-1}(B)\cap E;E,E\right).
\end{equation*}
Then in virtue of (\ref{eq5.7}) we
clearly have that
\begin{equation*}
\hat{f}_{\phi}(t,\psi(\tau))=\hat{h}(t,\tau)=\tilde{f}_{\psi}(\phi(t),\tau),\qquad
(t,\tau)\in\widehat{\X}\left(\phi^{-1}(A)\cap E,\psi^{-1}(B)\cap E;E,E\right).
\end{equation*}
Therefore,
\begin{equation*}
\hat{f}_{\phi}(0,w_0)=\hat{h}(0,0)=\tilde{f}_{\psi}(z_0,0),
\end{equation*}
which, in turn, implies that
\begin{equation*}
\hat{f}(z_0,w_0)=\hat{f}_{\phi}(0,w_0)=\tilde{f}_{\psi}(z_0,0)=\hat{g}(z_0,w_0).
\end{equation*}
Hence, the proof of identity (\ref{eq5.5}) is complete. This finishes  the proof of  the theorem.
\end{proof}

\section{Part 3 of the proof:  Theorem A for the case $N=2$ }

The main purpose of the section is  to prove Theorem A for the case $N=2.$
\begin{thm}\label{thmA_special_case}
Let  $D,$ $G$ be  complex manifolds,    and let $ A\subset
D,$  $B\subset G$  be non locally pluripolar  subsets.  Let $Z$ be a complex analytic space
possessing the Hartogs extension property. Put $X:=\X(A,B;D,G)$ and
$\widehat{X}:=\widehat{\X}(A,B;D,G).$
Then for any mapping $f\in\mathcal{O}_s(X,Z),$ there is a unique   mapping
$\hat{f}\in\mathcal{O}(\widehat{X},Z)$
such that $\hat{f}=f$ on $X\cap  X^{\ast}.$
\end{thm}

 For the proof we need to develop some preparatory results.
 
For any $a\in A^{\ast}$ (resp. $b\in B^{\ast}$), fix an
open neighborhood $U_{a}$ of  $a$  (resp.  $V_b$ of $b$) such that $U_{a} $
(resp. $V_b$) is  biholomorphic to a domain in $\C^{d_a}$  (resp. in
$\C^{d_b}$),
where $d_a $ (resp.   $d_b $) is the dimension of $D$ (resp. $G$) at $a$ (resp. $b$).
 For any $0<\delta\leq\frac{1}{2},$  define
\begin{equation}\label{eq6.1}
\begin{split}
U_{a,\delta}&:=\left\lbrace z\in U_{a}:\ \widetilde{\omega}(z, A\cap U_a,   U_a)<\delta  \right\rbrace,\qquad
a\in A\cap A^{\ast},\\
V_{b,\delta}&:=\left\lbrace w\in V_{b}:\  \widetilde{\omega}(w, B\cap V_b,   V_b)<\delta  \right\rbrace,\qquad
b\in B\cap B^{\ast},\\
A_{\delta}&:=\bigcup\limits_{a\in A\cap A^{\ast}} U_{a,\delta},\qquad
B_{\delta}:=\bigcup\limits_{b\in B\cap B^{\ast}} V_{b,\delta},\\
 D_{\delta}&:=\left\lbrace z\in D:\
  \widetilde{\omega}(z,A,D)<1-\delta\right\rbrace,\quad
  G_{\delta}:=\left\lbrace w\in G:\  \widetilde{\omega}(w,B,G)<1-\delta\right\rbrace.
\end{split}
\end{equation}

\begin{lem}\label{lem6.2}
 We keep the above notation.
Then
 \begin{eqnarray}
\label{eq6.2} A\cap A^{\ast}&\subset& A_{\delta}\subset D_{1-\delta} \subset D_{\delta},\\
\label{eq6.3} \widetilde{\omega}(z,A,D)-\delta
&\leq&  \widetilde{\omega}(z, A_{\delta},D) \leq\widetilde{\omega}(z,A,D),\qquad z\in D .
\end{eqnarray}
\end{lem}

\smallskip

\noindent {\it Proof of Lemma \ref{lem6.2}.} Using (\ref{eq6.1}) and  the definition
of local pluriregularity, we see that $a\in U_{a,\delta}$
for $a\in A\cap A^{\ast}.$ Consequently, the
first inclusion in (\ref{eq6.2})  follows. Since $0<\delta\leq\frac{1}{2},$ the third inclusion in (\ref{eq6.2}) is clear. To prove the second
inclusion in  (\ref{eq6.2}), let $z$ be an arbitrary point of
$A_{\delta}.$ Then there is an $a\in A\cap  A^{\ast}$ such that $z\in
U_{a,\delta}.$ Applying Part 2) of Proposition \ref{prop3.5}
and taking into account the inequality $0<\delta\leq\frac{1}{2},$  we
 obtain
\begin{equation}\label{eq6.4}
 \widetilde{\omega}(z,A,D)     
\leq  \widetilde{\omega}(z,A\cap U_a,U_a)<\delta .
\end{equation}
Hence, $z\in D_{1-\delta},$ which in turn implies that $A_{\delta}\subset
D_{1-\delta}.$   Hence, all assertions in (\ref{eq6.2}) are proved.

Next,  using the first inclusion in
(\ref{eq6.2}) and applying Parts 1) and 2) of Proposition \ref{prop3.5}, we get
 \begin{equation*}
 \widetilde{\omega}(z,A_{\delta},D)
\leq  \widetilde{\omega}(z, A\cap A^{\ast},D) = \widetilde{\omega}(z,A,D),\qquad z\in D,
\end{equation*}
which proves the second estimate in  (\ref{eq6.3}).

To complete the proof of  (\ref{eq6.3}), let   $a\in A\cap A^{\ast}$  and
$0<\delta\leq\frac{1}{2}.$ We deduce from (\ref{eq6.4})  that $
 \widetilde{\omega}(z,A ,D)-\delta\leq 0$ for  $   z\in U_{a,\delta}.$ Hence, by
(\ref{eq6.1}),
\begin{equation*}
 \widetilde{\omega}(z,A ,D)-\delta\leq 0, \qquad z\in
A_{\delta}.
\end{equation*}
On the other hand,  $ \widetilde{\omega}(z,A ,D)-\delta<1,$ $z\in D.$
Consequently, the first estimate in   (\ref{eq6.3})
follows. Hence, the proof of the  lemma is   finished.
\hfill $\square$

We also need the following
\begin{defi}\label{defi6.3}
Let $\mathcal{M}$ be a complex manifold and  $Y$  a complex space. Let $(U_j)_{j\in J}$ be
a family of open subsets of $\mathcal{M},$ and 
 $(f_j)_{j\in J}$ 
a family of   mappings such that $f_j\in \mathcal{O}(U_j,Y).$
We say that the family $(f_j)_{j\in J}$ is "collective" if, for any $j,k\in J,$    $f_j=f_k$ on $U_j\cap U_k.$
The unique holomorphic mapping $f: \bigcup\limits_{j\in J}U_j\longrightarrow Y,$ defined by
$f:=f_j$ on $U_j,$  $j\in J,$ is called the "collected mapping" of $(f_j)_{j\in J}.$
\end{defi}

\begin{lem}  \label{collecting_lem}
We keep the  the hypothesis of  Theorem \ref{thmA_special_case} and the above notation.
Suppose, in addition,  that  for every $a\in A\cap A^{\ast},$ there is a (unique) mapping
   $\hat{f}_{a} \in   \mathcal{O}\Big(\widehat{\X}
 \left(A\cap U_{a} ,B; U_{a},G\right),Z\Big) $
such that
\begin{equation}\label{eq6.5}
\hat{f}_{a}(z,w)= f(z,w),\  (z,w)\in
 \X\left(  A\cap A^{\ast}\cap  U_a,B\cap B^{\ast};  U_{a},G\right) .
\end{equation}
Then   the family $\left(\hat{f}_{a}|_{U_{a,\delta}\times G_{\delta}} \right)_{a\in A\cap A^{\ast}}$
 is collective.
\end{lem}

\smallskip

\noindent {\it Proof of Lemma \ref{collecting_lem}.} Let $a_1,\ a_2 $ be arbitrary  elements of $ A\cap A^{\ast}.$
  By (\ref{eq6.5}), we have   that  
\begin{equation*} 
\hat{f}_{a_1} (z,w)=f(z,w)=\hat{f}_{a_2}(z,w),\qquad (z,w)\in (U_{a_1}\cap U_{a_2})\times (B\cap B^{\ast}).
 \end{equation*}
Consequently, in virtue of Part 1) of Theorem  \ref{unique2},
    \begin{equation*}
\hat{f}_{a_1} (z,w)=\hat{f}_{a_2}(z,w),\qquad (z,w)\in\widehat{\X} \left(A\cap U_{a_1} ,B; U_{a_1},G\right)
\cap\widehat{\X}\left(A\cap U_{a_2} ,B; U_{a_2},G\right).
 \end{equation*}
  This,   combined with 
  the definition of $U_{a,\delta}$ and $G_{\delta}$  given in  (\ref{eq6.1}), the fact that $0<\delta\leq\frac{1}{2},$ and
  Definition \ref{defi6.3}, implies 
 the desired conclusion. 
\hfill $\square$

\begin{lem}  \label{limit_lem}
 Let $\mathcal{D}$ and $\mathcal{G}$ be two complex manifolds.
Let $(\mathcal{A}_{\delta})_{0<\delta<\frac{1}{2}}$  (resp. $(\mathcal{B}_{\delta})_{0<\delta<\frac{1}{2}}$)
be a family of non locally pluripolar subsets of $\mathcal{D}$ (resp. $\mathcal{G}$),
and  $(\mathcal{D}_{\delta})_{0<\delta<\frac{1}{2}}$  (resp. $(\mathcal{G}_{\delta})_{0<\delta<\frac{1}{2}}$)
  a family of  open subsets of $\mathcal{D}$ (resp. $\mathcal{G}$)  with the following properties:
\begin{itemize}
\item[(i)]  $\mathcal{A}_{\delta_1}\subset \mathcal{A}_{\delta_2}\subset \mathcal{D}_{\delta_2}\subset \mathcal{D}_{\delta_1}$  and 
$\mathcal{B}_{\delta_1}\subset \mathcal{B}_{\delta_2}
\subset \mathcal{G}_{\delta_2}\subset \mathcal{G}_{\delta_1}$ for $0<\delta_1\leq\delta_2<\frac{1}{2}.$
\item[(ii)]   There is a family of holomorphic mappings  $(\hat{f}_{\delta})_{0<\delta<\frac{1}{2}}$       such that
   $\hat{f}_{\delta} \in   \mathcal{O}\Big(\widehat{\X}
 \left(\mathcal{A}_{\delta} ,  \mathcal{B}_{\delta}  ;   \mathcal{D}_{\delta}, \mathcal{G}_{\delta}\right),Z\Big), $ and
 for $0<\delta_1<\delta_2<\frac{1}{2},$
 \begin{equation*} 
\hat{f}_{\delta_1}(z,w)= \hat{f}_{\delta_2}(z,w) ,\qquad  (z,w)\in
   \mathcal{A}_{\delta_1}\times \mathcal{B}_{\delta_1}.
\end{equation*}
\item[(iii)] 
 There are an open subset  $U$ (resp. $V$)  of $\mathcal{D}$  (resp. $\mathcal{G}$) and a number $0<\delta_0<\frac{1}{2}$  such that 
  $
 \widetilde{\omega}(z, \mathcal{A}_{\delta},\mathcal{D}_{\delta_0})+
 \widetilde{\omega}(w, \mathcal{B}_{\delta},\mathcal{G}_{\delta_0})<1
 $
for  all $(z,w)\in U\times V$ and $0<\delta <\delta_0.$
\end{itemize}
 Then   $
\hat{f}_{\delta}(z,w)= \hat{f}_{\delta_0}(z,w) $ for all $  (z,w)\in
    U\times  V $ and  $ 0<\delta<\delta_0.$
\end{lem}

\smallskip

\noindent {\it Proof of Lemma \ref{limit_lem}.}
 Fix  $\delta $ such that $ 0<\delta <\delta_0.$
By (iii), we have that 
\begin{equation}\label{eq_1_limit_lem}
U\times V \subset H:= \widehat{\X}\left(\mathcal{A}_{\delta} ,\mathcal{B}_{\delta};\mathcal{D}_{\delta_0},\mathcal{G}_{\delta_0}\right).
\end{equation}
 On the other hand, using (i) and  Part 2) of Proposition \ref{prop3.5},  we see that
\begin{equation*} 
 H \subset   \widehat{\X}\left(\mathcal{A}_{\delta} ,\mathcal{B}_{\delta};\mathcal{D}_{\delta},\mathcal{G}_{\delta}\right)\cap
   \widehat{\X}\left(\mathcal{A}_{\delta_0} ,\mathcal{B}_{\delta_0};\mathcal{D}_{\delta_0},\mathcal{G}_{\delta_0}\right).
\end{equation*}
 Using this and     (ii),  we are able to apply
   Part 2) of Theorem  \ref{unique2} to
 $\hat{f}_{\delta}|_H$ and $\hat{f}_{\delta_0}|_H.$
 Consequently,
 $\hat{f}_{\delta} =\hat{f}_{\delta_0} $ on $H.$  
This, combined with (\ref{eq_1_limit_lem}),  completes the proof of the lemma. \hfill
$\square$
 
Now we are able to  to prove Theorem  \ref{thmA_special_case} in the following
special case.
\begin{prop}\label{prop_special_case}
 We keep the   hypothesis of  Theorem \ref{thmA_special_case}. Suppose in addition that $G$ is biholomorphic
to a domain in $\C^q$ $(q\in\N).$  Then   the conclusion of  Theorem \ref{thmA_special_case} holds.
\end{prop}

\smallskip

\noindent {\it Proof of  Proposition \ref{prop_special_case}.}
   For each $a\in A\cap A^{\ast},$  let $f_a:=f|_{\X\left( A\cap U_{a} ,B;  U_{a},G\right)}.$
Since $f\in \mathcal{O}_s(X,Z),$ we deduce that    $f_a\in
\mathcal{O}_s\Big(\X\left( A\cap U_{a} ,B;  U_{a},G\right),Z\Big).$ Recall that
$U_a$ (resp. $G$) is biholomorphic to a domain in $\C^{d_a}$ (resp. in $\C^q$). Consequently,  applying Theorem
3
to $f_a$ yields that there is a unique mapping
   $\hat{f}_{a} \in   \mathcal{O}\Big(\widehat{\X}
 \left(A\cap U_{a} ,B; U_{a},G\right),Z\Big) $
such that
\begin{equation}\label{eq1_prop_special_case}
\hat{f}_{a}(z,w)=f_a(z,w)=f(z,w),\  (z,w)\in
 \X\left(  A\cap A^{\ast}\cap  U_a,B\cap B^{\ast};  U_{a},G\right) .
\end{equation}
Let $0<\delta\leq \frac{1}{2}.$ 
In virtue of  (\ref{eq1_prop_special_case}), we are able to apply Lemma  \ref{collecting_lem} to the family $\left(\hat{f}_{a}|_{U_{a,\delta}\times G_{\delta}} \right)_{a\in A\cap A^{\ast}}.$ 
Let
\begin{equation}\label{eq2_prop_special_case}
\tilde{\tilde{f}}_{\delta}\in \mathcal{O}(A_{\delta}\times
G_{\delta},Z)
\end{equation}
denote the collected mapping of this family.
 In virtue of  (\ref{eq1_prop_special_case})--(\ref{eq2_prop_special_case}),
 we are able to define a new mapping  $\tilde{f}_{\delta}$ on $\X\left(A_{\delta},B\cap B^{\ast};D,G_{\delta}
 \right)$ as follows
\begin{equation*}
 \tilde{f}_{\delta}:=
\begin{cases}
 \tilde{\tilde{f}}_{\delta},
  & \qquad\text{on}\  A_{\delta}\times G_{\delta}, \\
  f, &   \qquad\text{on}\ D\times (B\cap B^{\ast})        .
\end{cases}
\end{equation*}
Using this  and  (\ref{eq1_prop_special_case})--(\ref{eq2_prop_special_case}), we  see that  $ \tilde{f}_{\delta}\in \mathcal{O}_s\Big(
\X\left(A_{\delta}, B\cap B^{\ast};D,G_{\delta}
 \right),Z\Big),$ and
\begin{equation}\label{eq3_prop_special_case}
  \tilde{f}_{\delta}=f\qquad\text{on}\ \X(A\cap A^{\ast},B\cap
 B^{\ast};D,G_{\delta}).
\end{equation}
 Since $A_{\delta}$ is open and $G_{\delta}$ is biholomorphic to an open set in $\C^q,$  we are able to  apply  Theorem  \ref{thm4.1}
to  $ \tilde{f}_{\delta}$ in order to obtain a mapping
  $\hat{f}_{\delta}\in \mathcal{O}\Big(  \widehat{\X}\left(A_{\delta}, B\cap B^{\ast};D,G_{\delta}
 \right),Z \Big)$  such that
\begin{equation}\label{eq4_prop_special_case }
  \hat{f}_{\delta}=  \tilde{f}_{\delta}\qquad\text{on}\ \X\left(A_{\delta},
  B\cap B^{\ast};D,G_{\delta}\right).
\end{equation}

We are now in a position to define the desired extension  mapping $\hat{f}.$
Indeed, one  glues
$\left(\hat{f}_{\delta}\right)_{0<\delta\leq\frac{1}{2}}$ together to obtain
$\hat{f}$ in the following way
\begin{equation}\label{eq5_prop_special_case}
\hat{f}:=\lim\limits_{\delta\to 0} \hat{f}_{\delta}\qquad \text{on}\
 \widehat{X}  .
\end{equation}
One needs to check that the limit (\ref{eq5_prop_special_case}) exists and possesses all the required
properties.  In virtue of (\ref{eq3_prop_special_case})--(\ref{eq5_prop_special_case}),  the fact that
$G_{\delta}\nearrow G$ as $\delta\searrow 0$ (by (\ref{eq6.1})) and Lemma \ref{limit_lem},
the proof will be complete if we can show that for every
$(z_0,w_0)\in\widehat{X},$ there are an open neighborhood
$U\times V$ of  $(z_0,w_0)$ and   $\delta_0>0$ such that the hypothesis of Lemma \ref{limit_lem}
is fulfilled with 
\begin{equation*}
\mathcal{D}:=D,\ \mathcal{G}:=G,\ \mathcal{A}_{\delta}:=A_{\delta},\  \mathcal{B}_{\delta}:=B\cap B^{\ast},\
\mathcal{D}_{\delta}:=D,\ \mathcal{G}_{\delta}:=G_{\delta},\qquad 0<\delta<\frac{1}{2}.
\end{equation*}
  
 To this end let  
 \begin{equation} \label{eq6_prop_special_case}
 \delta_{0}:=\frac{1-\widetilde{\omega}(z_0,A,D)- \widetilde{\omega}(w_0,B,G)}{2},
\end{equation}
 and let $U\times V$ be an open  neighborhood of  $(z_0,w_0)$ such that
 \begin{equation} \label{eq7_prop_special_case}
   \widetilde{\omega}(z,A,D)+ \widetilde{\omega}(w,B,G)< \widetilde{\omega}(z_0,A,D)+ \widetilde{\omega}(w_0,B,G)+\delta_0.
\end{equation}
Then for $0<\delta <\delta_0$ and for $(z,w)\in U\times V,$    using  (\ref{eq6_prop_special_case})--(\ref{eq7_prop_special_case}) and invoking Part 4) of Proposition \ref{prop3.5}, we see that
\begin{equation}
\begin{split}\label{eq8_prop_special_case}
 \widetilde{\omega}(z,A_{\delta},D)+ \widetilde{\omega}(w,B\cap B^{\ast},G_{\delta_0})&\leq
 \widetilde{\omega}(z,A,D)+\frac{ \widetilde{\omega}(w,B,G)}{1-\delta_0}\\
&\leq\frac{ \widetilde{\omega}(z,A,D)+ \widetilde{\omega}(w,B,G) }{1-\delta_0}
<1.
\end{split}
\end{equation}
This proves the above assertion. Hence, the proof of the proposition is finished.
\hfill  $\square$

\smallskip

We now arrive at

\smallskip

\noindent {\it Proof of  Theorem \ref{thmA_special_case}.}
  For each $a\in A\cap A^{\ast},$  let $f_a:=f|_{\X\left( A\cap U_{a} ,B;  U_{a},G\right)}.$
Since $f\in \mathcal{O}_s(X,Z),$ we deduce that    $f_a\in
\mathcal{O}_s\Big(\X\left( A\cap U_{a} ,B;  U_{a},G\right),Z\Big).$ Since $U_a$ is
biholomorphic to a domain in $\C^{d_a},$ we are able to apply  Proposition \ref{prop_special_case}
 to $f_a.$  Consequently,  there is a unique mapping
   $\hat{f}_{a} \in   \mathcal{O}\Big(\widehat{\X}
 \left(A\cap U_{a} ,B; U_{a},G\right),Z\Big) $
such that
\begin{equation}\label{eq1_thmA_special_case}
\hat{f}_{a}(z,w)=f(z,w),\qquad  (z,w)\in
 \X\left(  A\cap A^{\ast}\cap U_a,B\cap B^{\ast};  U_{a},G\right) .
\end{equation}
 Let $0<\delta\leq\frac{1}{2}.$ In virtue of  (\ref{eq1_thmA_special_case}), we may apply Lemma  \ref{collecting_lem}.   Consequently, we can collect  the family
 $\left(\hat{f}_{a}|_{U_{a,\delta}\times G_{\delta}} \right)_{a\in A\cap A^{\ast}}$
in order to obtain the collected   mapping $\tilde{f}^A_{\delta}\in \mathcal{O}(A_{\delta}\times
G_{\delta},Z).$

Similarly,  for each $b\in B\cap B^{\ast},$   one obtains a unique mapping
 $\hat{f}_{b} \in   \mathcal{O}\Big(\widehat{\X}
 \left(A, B\cap V_{b} ; D,V_{b}\right),Z\Big) $
such that
\begin{equation}\label{eq2_thmA_special_case}
\hat{f}_{b}(z,w) =f(z,w),\qquad (z,w)\in
 \X\left(  A\cap A^{\ast},B\cap B^{\ast}\cap V_b;  D,V_{b}\right) .
\end{equation}
Moreover, one can collect  the family
 $\left(\hat{f}_{b}|_{D_{\delta}\times V_{b,\delta}} \right)_{b\in B\cap B^{\ast}}$
in order to obtain the collected   mapping $\tilde{f}^B_{\delta}\in \mathcal{O}(D_{\delta}\times
B_{\delta},Z).$

Next, we prove that
\begin{equation} \label{eq3_thmA_special_case}
 \tilde{f}^A_{\delta}=\tilde{f}^B_{\delta}\qquad \text{on}\  A_{\delta}\times
B_{\delta}.
\end{equation}
Indeed, in virtue of  (\ref{eq1_thmA_special_case})--(\ref{eq2_thmA_special_case})  it suffices to show that
 for any $a\in A\cap A^{\ast}$ and $b\in B\cap B^{\ast}$ and any $0<\delta\leq\frac{1}{2},$  
\begin{equation}\label{eq4_thmA_special_case}
\hat{f}_{a} (z,w)=\hat{f}_{b}(z,w),\qquad (z,w)\in U_{a,\delta}\times V_{b,\delta}.
 \end{equation}
 Observe that in virtue of  (\ref{eq1_thmA_special_case})--(\ref{eq2_thmA_special_case})  
 one has that
\begin{equation*} 
 \hat{f}_{a}(z,w)=\hat{f}_{b}(z,w)=f(z,w), \qquad (z,w)\in   \X\left(A\cap A^{\ast}\cap U_{a}
,B\cap B^{\ast}\cap V_{b}; U_{a},V_b\right).
\end{equation*}
  Recall that $U_a$ (resp. $V_b$) is biholomorphic to a domain
in $\C^{d_a}$ (resp. $\C^{d_b}$). Consequently, applying the uniqueness of Theorem 3 yields that
 \begin{equation*}
\hat{f}_{a}(z,w)= \hat{f}_{b}(z,w) , \qquad (z,w)\in \widehat{\X}\left(A\cap U_{a} ,B\cap V_{b};
U_{a},V_b\right).
\end{equation*}
   Hence, the proof of  (\ref{eq4_thmA_special_case}) and then the proof of (\ref{eq3_thmA_special_case})          are finished.

 In virtue of  (\ref{eq3_thmA_special_case}),
 we are able to define a new  mapping $\tilde{f}_{\delta}:\
 \X\left(A_{\delta}, B_{\delta};D_{\delta},
 G_{\delta}\right)\longrightarrow Z$ as follows
\begin{equation}\label{eq5_thmA_special_case}  
 \tilde{f}_{\delta}:=
\begin{cases}
\tilde{f}^A_{\delta},
  & \qquad\text{on}\  A_{\delta}\times G_{\delta}, \\
  \tilde{f}^B_{\delta}, &   \qquad\text{on}\ D_{\delta}\times B_{\delta}        .
\end{cases}
\end{equation}
Using formula (\ref{eq5_thmA_special_case}) it can be readily checked that
$\tilde{f}_{\delta}\in \mathcal{O}_s\Big(\X\left(A_{\delta}, B_{\delta};D_{\delta},G_{\delta}
 \right),Z\Big).$ Since  we know from  (\ref{eq6.2}) that  $A_{\delta}$
 (resp.  $B_{\delta}$) is an open subset of $D_{\delta}$  (resp.
 $G_{\delta}$), we are able to apply  Theorem  \ref{thm5.1}
to  $ \tilde{f}_{\delta}$  for every $0<\delta\leq\frac{1}{2}.$ Consequently,
one
  obtains a unique  mapping $\hat{f}_{\delta}\in \mathcal{O}\Big(
 \widehat{\X}\left(A_{\delta}, B_{\delta};D_{\delta},G_{\delta}
 \right),Z \Big)$ such that
\begin{equation}\label{eq6_thmA_special_case}  
  \hat{f}_{\delta}=  \tilde{f}_{\delta}\qquad\text{on}\ \X\left(A_{\delta},
  B_{\delta};D_{\delta},G_{\delta}\right).
\end{equation}
It follows from   (\ref{eq1_thmA_special_case})--(\ref{eq2_thmA_special_case}) and   (\ref{eq5_thmA_special_case})--(\ref{eq6_thmA_special_case})  that
 \begin{equation}\label{eq7_thmA_special_case}
  \hat{f}_{\delta}=f\qquad\text{on}\
\X\left(A \cap A^{\ast}, B\cap B^{\ast};D_{\delta},G_{\delta}\right).
 \end{equation}
 In addition, for any $0<\delta\leq\delta_0\leq\frac{1}{2},$ and any $(z,w)\in
 A_{\delta}\times B_{\delta},$ there is an  $a\in A\cap A^{\ast}$ 
  such that $z\in U_{a,\delta_0}.$ 
  Therefore, it
 follows from the construction of $\tilde{f}^A_{\delta},$
    (\ref{eq5_thmA_special_case}) and (\ref{eq6_thmA_special_case}) that
\begin{equation*}
  \hat{f}_{\delta}(z,w)=  \hat{f}_a(z,w)=\hat{f}_{\delta_0}(z,w) .
\end{equation*}
This proves that
\begin{equation}\label{eq8_thmA_special_case}
  \hat{f}_{\delta} =  \hat{f}_{\delta_0} \qquad\text{on}\
  A_{\delta}\times B_{\delta},\ 0<\delta\leq\delta_0\leq\frac{1}{2}
   .
\end{equation}

We are now in a position to define the desired extension  mapping $\hat{f}.$
 \begin{equation*}
\hat{f}:=\lim\limits_{\delta\to 0} \hat{f}_{\delta}\qquad \text{on}\
 \widehat{\X}  .
\end{equation*}
 To prove that $\hat{f}$  satisfies   the desired conclusion of the theorem one proceeds as in the end of the proof of Proposition \ref{prop_special_case}.   In virtue of (\ref{eq7_thmA_special_case})--(\ref{eq8_thmA_special_case}) and Lemma \ref{limit_lem},
the proof will be complete if we can  verify that for every
$(z_0,w_0)\in\widehat{X},$ there are an open neighborhood
$U\times V$ of  $(z_0,w_0)$ and   $\delta_0>0$ such that the hypothesis of Lemma \ref{limit_lem}
is fulfilled with 
\begin{equation*}
\mathcal{D}:=D,\ \mathcal{G}:=G,\ \mathcal{A}_{\delta}:=A_{\delta},\  \mathcal{B}_{\delta}:=B_{\delta},\
\mathcal{D}_{\delta}:=D_{\delta},\ \mathcal{G}_{\delta}:=G_{\delta},\qquad 0<\delta<\frac{1}{2}.
\end{equation*}
  Since the verification follows along almost the same lines as  (\ref{eq6_prop_special_case})--(\ref{eq8_prop_special_case}),
  it is therefore left to the interested reader. 
 
Hence, the proof of the theorem  is finished.
 \hfill $\square$

\section{Part 4: Completion of the proof of Theorem A}

In this section we prove Theorem A for every $N\geq 3.$
  We divide the proof into two parts.

\subsection{Proof of the  existence and uniqueness of $\hat{f}$}

We proceed by induction (I) on $N\geq 3.$
Suppose the  theorem is true
for $N-1\geq 2.$ We have to discuss the case of an $N$-fold cross
$X:=\X(A_1,\ldots,A_N;D_1,\ldots,D_N),$ where $D_1,\ldots,D_N $ are
complex manifolds
  and $A_1\subset D_1,\ldots,A_N\subset D_N$ are   non locally pluripolar subsets
  $(1\leq j\leq N).$
Let  $f\in\mathcal{O}_s(X,Z).$ Observe that the uniqueness of $\hat{f}$ follows immediately
 from  Part 2) of Theorem \ref{unique2}.

We proceed again by induction (II) on the   integer $k$ $(0\leq k\leq N)$ such
that there are at least $k$ complex manifolds among
$\{D_1,\ldots,D_N\}$ which are  biholomorphic to  Euclidean domains.

For $k=N$ we are reduced to Theorem 3.

 Suppose that Theorem A is true
for the case where  $k=k_0$ $(1\leq k_0\leq N).$
 We have to discuss the case where $k=k_0-1.$ Suppose without loss of
 generality that $D_2$ is not biholomorphic to an Euclidean domain.

 For any $1\leq j\leq N$ and  $a_j\in A_j\cap A^{\ast}_j,$ one fixes an
open neighborhood $U_{a_j}$ of  $a_j$ such that $U_{a_j}$ is  biholomorphic to  a domain
in $\C^{d_{a_j}},$ where $d_{a_j}$ is the dimension of $D_j$ at $a_j.$
 For $1\leq j\leq N$ and for any $0<\delta< 1 ,$    define
\begin{equation}\label{eq8.1.1}
\begin{split}
U_{a_j,\delta}&:=\left\lbrace z_j\in U_{a_j}:\ \widetilde{\omega}(z_j, A_j\cap U_{a_j},  U_{a_j})<\delta
  \right\rbrace,\qquad
a_j\in A_j\cap A^{\ast}_j,\\
A_{j,\delta}&:=\bigcup\limits_{a_j\in A_j\cap A^{\ast}_j} U_{a_j,\delta},\\
D_{j,\delta}&:=\left\lbrace z_j\in D_j:\  \widetilde{\omega}(z_j, A_j ,   D_j)<1-\delta
  \right\rbrace.
\end{split}
\end{equation}

For every $a_1\in A_1\cap A^{\ast}_1,$ consider the mapping $f_{a_1}$
provided by
\begin{multline*}
f_{a_1}\left(z_2,\ldots,z_N\right):=f\left( a_1,z_2,\ldots,z_N \right),\\
\left(z_2,\ldots,z_N \right)\in\X\Big( A_2,\ldots,A_N;D_2,\ldots,D_N   \Big).
\end{multline*}
Observe that in virtue of   the above formula and the hypothesis that $f\in \mathcal{O}_s(X,Z),$
 $f_{a_1} $ satisfies
 the hypothesis of  Theorem A for $(N-1)$-cross.
 Consequently,   applying the hypothesis of induction (I), we obtain a
unique mapping $\hat{f}_{a_1} \in   \mathcal{O}\Big(\widehat{\X}
 \left(     A_2,\ldots,A_N; D_2,\ldots,D_N      \right),Z\Big) $
such that
\begin{multline}\label{eq8.1.2}
\hat{f}_{a_1}\left(z_2,\ldots,z_N\right)=f\left( a_1,z_2,\ldots,z_N \right),\\
\left(z_2,\ldots,z_N \right)\in\X\Big( A_2\cap A^{\ast}_2,\ldots,A_N\cap A^{\ast}_N
;D_2,\ldots,D_N   \Big).
\end{multline}

For every $a_2\in A_2\cap A^{\ast}_2,$ consider the mapping $f_{a_2}$
provided by
\begin{multline*}
f_{a_2}\left(z_1,z_2,z_3,\ldots,z_N\right):=f\left( z_1,z_2,z_3,\ldots,z_N \right),\\
\left(z_1,z_2,z_3,\ldots,z_N \right)\in\X\Big(A_1, A_2\cap U_{a_2},A_3,\ldots,A_N;D_1,U_{a_2},D_3,
\ldots,D_N   \Big).
\end{multline*}
Recall that  $U_{a_2}$ is biholomorphic to an Euclidean domain, but $D_2$
is not so. Therefore,   in virtue of   the above formula and the hypothesis that $f\in \mathcal{O}_s(X,Z),$
 we may apply the hypothesis of induction (II) to $f_{a_2}. $
 Consequently,   one obtains a
unique mapping $\hat{f}_{a_2} \in   \mathcal{O}\Big(\widehat{\X}
 \left(     A_1,A_2\cap U_{a_2},A_3,\ldots,A_N;D_1, U_{a_2}, D_3,\ldots,D_N      \right),Z\Big) $
such that
\begin{multline}\label{eq8.1.3}
\hat{f}_{a_2}\left(z_1,z_2,z_3,\ldots,z_N\right)=f\left( z_1,z_2,z_3,\ldots,z_N \right),\\
\left(z_1,z_2,z_3,\ldots,z_N \right)\in\X\Big( A_1\cap A^{\ast}_1,A_2\cap A^{\ast}_2\cap U_{a_2},
A_3\cap A^{\ast}_3,\ldots,A_N\cap A^{\ast}_N
;D_1, U_{a_2},D_3,\ldots,D_N   \Big).
\end{multline}

 We need the following
\begin{lem}\label{lem8.1}
We keep the hypothesis of Theorem A and the above notation.
Then for any $a_1\in A_1\cap A^{\ast}_1,$  $a_2\in A_2\cap A^{\ast}_2,$
and any $0<\delta < \frac{1}{N},$
\begin{multline*}
\hat{f}_{a_1}\left(z_2,z_3,\ldots,z_N\right)=\hat{f}_{a_2}\left( a_1,z_2,z_3,\ldots,z_N \right),\\
\left(z_2,z_3,\ldots,z_N \right)\in U_{a_2,\delta}\times A_{3,\delta}\times\cdots
\times A_{N,\delta}.
\end{multline*}
\end{lem}

\smallskip

\noindent {\it Proof of Lemma \ref{lem8.1}.} Let $a_1,$ $a_2$ be as in the
statement of Lemma \ref{lem8.1}. In virtue of (\ref{eq8.1.2})--(\ref{eq8.1.3}),
\begin{multline*}
\hat{f}_{a_2}\left(a_1,z_2,z_3,\ldots,z_N\right)=f\left( a_1,z_2,z_3,\ldots,z_N \right)=\hat{f}_{a_1}\left(z_2,z_3,\ldots,z_N\right),\\
\left(z_2,z_3,\ldots,z_N \right)\in\X\Big( A_2\cap A^{\ast}_2\cap U_{a_2},
A_3\cap A^{\ast}_3,\ldots,A_N\cap A^{\ast}_N
;  U_{a_2},D_3,\ldots,D_N   \Big).
\end{multline*}
 Consequently, applying  Part 2) of Theorem \ref{unique2}
to $\hat{f}_{a_1}$ and $\hat{f}_{a_2}(a_1,\cdot)$  yields that
 \begin{multline}\label{eq8.1.5}
\hat{f}_{a_1}\left(z_2,z_3,\ldots,z_N\right)=\hat{f}_{a_2}\left( a_1,z_2,z_3,\ldots,z_N \right),\\
\left(z_2,z_3,\ldots,z_N \right)\in
\widehat{\X}\Big( A_2 \cap U_{a_2},
A_3 ,\ldots,A_N
;  U_{a_2},D_3,\ldots,D_N   \Big).
\end{multline}
Moreover,  since $0<\delta<\frac{1}{N},$ it follows from (\ref{eq8.1.1}),  (\ref{eq6.2}), and a
straightforward computation that
\begin{multline*}
 U_{a_2,\delta}\times A_{3,\delta}\times\cdots
\times A_{N,\delta}\subset
\widehat{\X}\Big( A_2 \cap U_{a_2},
A_3 ,\ldots,A_N
;  U_{a_2},D_3,\ldots,D_N   \Big).
\end{multline*}
This, combined with (\ref{eq8.1.5}), implies the desired conclusion of
the lemma. \hfill $\square$

\smallskip

 In the sequel we always suppose that $0<\delta<\frac{1}{N}.$ In virtue of
 (\ref{eq8.1.3}),
    Part 1) of Theorem \ref{unique2} and Definition \ref{defi6.3}, we are able to
collect the family of mappings
\begin{equation*}
\left(\hat{f}_{a_2}|_{    D_{1,N\delta}\times U_{a_2,\delta}\times A_{3,\delta}\times \cdots
   \times A_{N,\delta}}
\right)_{a_2\in A_2\cap A^{\ast}_2}
\end{equation*}
 in order to obtain the collected  mapping
 \begin{equation}\label{eq8.1.6}
 \tilde{\tilde{f}}_{\delta}\in \mathcal{O}\Big( D_{1,N\delta}\times A_{2,\delta}\times\cdots
   \times A_{N,\delta} ,Z  \Big ).
   \end{equation}
    Let
    \begin{equation}\label{eq8.1.7}
    \begin{split}
 X_{\delta}&:= \X\Big(A_1\cap A^{\ast}_1, A_{2,\delta}\times\cdots \times  A_{N,\delta};D_{1,N\delta},
 \widehat{\X}\left(A_2,\ldots,A_N;D_2 ,\ldots, D_N\right)\Big),\\
\widehat{X}_{\delta}&:=\widehat{\X}\Big(A_1\cap A^{\ast}_1, A_{2,\delta}\times\cdots \times  A_{N,\delta};D_{1,N\delta},
     \widehat{ \X}\left(A_2,\ldots,A_N;D_2,\ldots,D_N\right) \Big).
 \end{split}
  \end{equation}
In virtue of Lemma \ref{lem8.1} and the construction (\ref{eq8.1.6}),
 we are able to define a new mapping $\tilde{f}_{\delta}:  \ X_{\delta}\longrightarrow
 Z$ as follows
\begin{equation}\label{eq8.1.8}
 \tilde{f}_{\delta}(z ):=
\begin{cases}
\tilde{\tilde{f}}_{\delta}(z),
  & \qquad z\in D_{1,N\delta}\times A_{2,\delta}\times\cdots   \times A_{N,\delta} ,\\
  \hat{f}_{z_1}, &
    \qquad z\in (A_1\cap A^{\ast}_1)\times\widehat{
    \X}(A_2,\ldots,A_N;D_2,\ldots,D_N),
\end{cases}
\end{equation}
where $z=(z_1,\ldots,z_N)\in X_{\delta}.$

Using (\ref{eq8.1.8}), (\ref{eq8.1.2}) and (\ref{eq8.1.6}),   it can be readily checked that
 $\tilde{f}_{\delta}\in \mathcal{O}_s(X_{\delta}).$ In addition,
 using (\ref{eq8.1.7}) we have that $ X_{\delta}\cap
X_{\delta}^{\ast}= X_{\delta}.$
Consequently, for every $0<\delta<\frac{1}{N},$  one applies Theorem \ref{thmA_special_case}
to  $\tilde{f}_{\delta}$
 and obtain
a  unique mapping  $\hat{f}_{\delta}\in
\mathcal{O}\left(\widehat{X}_{\delta}\right)$ such that
\begin{equation}\label{eq8.1.9}
\hat{f}_{\delta}=
\tilde{f}_{\delta} \qquad\text{on}\ X_{\delta}.
\end{equation}

 Finally, gluing $\left(\hat{f}_{\delta}\right)_{0<\delta<\frac{1}{N}},$ we can define the
   desired extension  mapping $\hat{f}$  by the formula
\begin{equation}\label{eq8.1.10}
\hat{f}:=\lim\limits_{\delta\to 0} \hat{f}_{\delta}\qquad \text{on}\
\widehat{X} .
\end{equation}

Next, we argue as in the proof of Theorem \ref{thmA_special_case}. More precisely, one  checks that the  hypothesis of Lemma \ref{limit_lem}
is fulfilled with  
 \begin{multline*}
 \mathcal{D}:=D_1,\ \mathcal{G} :=\widehat{ \X}\left(A_2,\ldots,A_N;D_2,\ldots,D_N\right),\\
 \mathcal{A}_{\delta}:=A_1\cap A_1^{\ast},\  
 \mathcal{B}_{\delta}:=A_{2,\delta}\times\cdots \times  A_{N,\delta},\ \mathcal{D}_{\delta}:=D_{1,N\delta},\  \mathcal{G}_{\delta}:=\widehat{ \X}\left(A_2,\ldots,A_N;D_2,\ldots,D_N\right), 
 \end{multline*}
for  $0<\delta<\frac{1}{N}.$

To do this let
\begin{eqnarray*}
\Omega&:=&\widehat{\X}\left(A_{2,\delta},\ldots,A_{N,\delta};D_2,\ldots,D_N\right),\\
\Omega_{N\delta}&:=&\left\lbrace z^{'}\in \Omega:\ \widetilde{\omega}\Big( z^{'},
A_{2,\delta}\times\cdots \times A_{N,\delta},\Omega\Big) <1-N\delta  \right\rbrace.
\end{eqnarray*}
Applying inequality (\ref{eq6.3}), one gets
\begin{equation*}
\sum\limits_{j=2}^N\widetilde{\omega}(z_j,A_{j},D_j)<\sum\limits_{j=2}^N\widetilde{\omega}(z_j,A_{j,\delta},D_j)+N\delta <1,\qquad
z^{'}=(z_2,\ldots,z_N)\in \Omega_{N\delta}.
\end{equation*}
Hence,  $\Omega_{N\delta}\subset
\widehat{\X}\left(A_2,\ldots,A_N;D_2,\ldots,D_N\right),$ which, in virtue of Part 2) of Proposition
\ref{prop3.5}, implies that
\begin{multline}\label{eq8.1.14}
\widetilde{\omega}\Big( z^{'},
A_{2,\delta}\times\cdots \times A_{N,\delta},\widehat{\X}\left(A_2,\ldots,A_N;D_2,\ldots,D_N\right)
       \Big) \\
       \leq \widetilde{\omega}\Big( z^{'},
A_{2,\delta}\times\cdots \times A_{N,\delta},\Omega_{N\delta}\Big),\qquad z^{'}\in \Omega_{N\delta}.
\end{multline}
On the other hand, in virtue of Part 2) of Proposition \ref{prop3.7}, we have that
\begin{equation}\label{eq8.1.15} \widetilde{\omega}\Big( z^{'},
A_{2,\delta}\times\cdots \times A_{N,\delta},\Omega\Big) =
\sum\limits_{j=2}^N\widetilde{\omega}(z_j,A_{j,\delta},D_j) ,\qquad z^{'}=(z_2,\ldots,z_N)\in \Omega.
\end{equation}
 By Part 4) of Proposition \ref{prop3.5}, 
\begin{equation*}  \widetilde{\omega}\Big( z^{'},
A_{2,\delta}\times\cdots \times A_{N,\delta},\Omega_{N\delta}\Big) =
 \frac{ \widetilde{\omega}\Big( z^{'},
A_{2,\delta}\times\cdots \times A_{N,\delta},\Omega\Big) }{1-N\delta},\qquad z^{'}\in \Omega_{N\delta}.
\end{equation*}
This, combined with (\ref{eq8.1.14})--(\ref{eq8.1.15}),  implies that
\begin{equation}\label{eq8.1.16}
\widetilde{\omega}\Big(z^{'},A_{2,\delta}\times\cdots \times A_{N,\delta},
 \widehat{\X}\left(A_2,\ldots,A_N;D_2,\ldots,D_N\right) \Big) 
 \leq \frac{ \sum\limits_{j=2}^N\widetilde{\omega}\Big( z_j,
A_{j,\delta} , D_j\Big) }{1-N\delta}.
 \end{equation}
 
 For every  $z^0=(z^0_1,z^{0'})\in\widehat{X},$  let
  $\delta_{0}:=  \frac{ 1-\sum\limits_{j=1}^N\widetilde{\omega}( z^0_j,
A_{j} , D_j) }{ 2N }     $  and fix an open neighborhood $U\times V$ of $z^0$ such that
\begin{equation*} \sum\limits_{j=1}^N\widetilde{\omega}( z_j,
A_{j} , D_j)<\delta_0+ \sum\limits_{j=1}^N\widetilde{\omega}( z^0_j,
A_{j} , D_j),\qquad z=(z_1,z^{'})\in U\times V.
\end{equation*}
 Then, using the latter estimate and  (\ref{eq8.1.16}) and Part 4) of Proposition \ref{prop3.5}, we see that 
\begin{multline*}
\widetilde{\omega}\left(z_1,A_1\cap A^{\ast}_1,D_{1,N\delta_0}\right)+\widetilde{\omega}\Big(z^{'},A_{2,\delta}\times\cdots \times A_{N,\delta},
 \widehat{\X}\left(A_2,\ldots,A_N;D_2,\ldots,D_N\right) \Big)\\
 \leq
  \frac{ \sum\limits_{j=1}^N\widetilde{\omega}( z^0_j,
A_{j} , D_j)+\delta_0 }{1-N\delta_0}<1 .
\end{multline*}
for  $z=(z_1,z^{'})\in U\times V$ and $0<\delta\leq\delta_0 .$
 Consequently, we are able to apply Lemma \ref{limit_lem}.

We complete the proof as follows.
An immediate consequence of Lemma \ref{limit_lem} and formula (\ref{eq8.1.10}) is that
$\hat{f}\in\mathcal{O}\left(\widehat{X},Z  \right).$
  Moreover,  by  (\ref{eq8.1.2})--(\ref{eq8.1.3}) and
  (\ref{eq8.1.7})--(\ref{eq8.1.10}), and using  the fact that
  $D_{1,N\delta}\nearrow D_1$ as $\delta\searrow 0$  (see (\ref{eq8.1.1})), we conclude   that $\hat{f}=f$  on the following set
   \begin{equation*}
 (A_1\cap A^{\ast}_1)\times \X\left(A_2\cap A^{\ast}_2,\ldots,
   A_N\cap A^{\ast}_N; D_2,\ldots,D_N  \right) \bigcup
  D_1\times  (A_2\cap A^{\ast}_2)\times  \cdots\times
   (A_N\cap A^{\ast}_N ) .
\end{equation*}
 Since this set  is equal to $X\cap
 X^{\ast}, $ it follows from Theorem \ref{unique3} that  the mapping $\hat{f}$ provided by formula (\ref{eq8.1.10})
 possesses all the desired properties. This completes induction (II) for
 $k=k_0-1.$
Hence, the proofs of induction (II), induction (I) and then the first part of the theorem  are finished. \hfill
$\square$
\subsection{Proof of the estimate in Theorem A}
Following the work in \cite{pn1} we divide this part into two steps.

\smallskip

\noindent {\bf Step 1:}
{\it Proof of the inequality $\vert \hat{f}\vert_{\widehat{X}}\leq\vert f\vert_X.$}

\smallskip

\noindent {\it Proof of Step 1.}
In order to reach a contradiction assume   that there is a point
$z^0\in\widehat{X}$ such that $\vert\hat{f}(z^0)\vert>\vert
f\vert_X.$ Put $\alpha:=\hat{f}(z^0)$ and consider the function
\begin{equation}\label{eq8.2.1}
g(z):=\frac{1}{f(z)-\alpha},\qquad  z\in X.
\end{equation}
Using the above assumption, we clearly have that  $g\in\mathcal{O}_s(X,\C).$ Hence by
 Subsection 8.1, there is exactly
one function  $\hat{g}\in\mathcal{O}(\widehat{X},\C)$ with $\hat{g}=g$ on
$X.$ Therefore, by (\ref{eq8.2.1}) we have on $X:$ $g(f-\alpha)\equiv 1.$ Thus
$\hat{g}(\hat{f}-\alpha)\equiv 1$ on $\widehat{X}.$ In particular,
\begin{equation*}
0=\hat{g}(z^0 )(\hat{f}(z^0)-\alpha)= 1,
\end{equation*}
which is a contradiction.
Hence the inequality $\vert \hat{f}\vert_{\widehat{X}}\leq\vert f\vert_X$ is
proved.    Thus Step 1 is complete. \hfill $\square$

\smallskip

\noindent {\bf Step 2:}
{\it Proof of the inequality
\begin{equation}\label{eq8.2.2}
\vert \hat{f}(z)\vert\leq\vert f\vert_A^{1-\omega(z)}\vert f\vert_X^{\omega(z)}
.
\end{equation}}
\noindent {\it Proof of Step 2.} We prove (\ref{eq8.2.2}) by induction on $N.$
 When $N=1,$ then  applying  Theorem
 \ref{two-constant} to the plurisubharmonic function $z\in D_1\mapsto \log\vert \hat{f}(z)\vert,$
 (\ref{eq8.2.2}) follows.
Suppose that (\ref{eq8.2.2})  is true for $N-1.$ We would like to prove it for
$N.$ Fix an arbitrary point
$z^0=\left(z^0_1,\ldots,z^0_N\right)\in\widehat{X}.$ Let
\begin{equation}\label{eq8.2.3}
\delta:=\sum\limits_{j=2}^N\widetilde{\omega}(z^0_j,A_j,D_j).
\end{equation}

For any $a_1\in A_1\cap A^{\ast}_1,$ we apply the hypothesis of induction
to the function $\hat{f}_{a_1}$ and obtain the estimate
\begin{equation}\label{eq8.2.4}
\vert \hat{f}_{a_1}\left(z^0_2,\ldots,z^0_N\right)\vert\leq\vert f\vert_A^{1-\delta}\vert f\vert_X^{\delta}
.
\end{equation}
 In virtue of (\ref{eq8.1.8})--(\ref{eq8.1.10}), we  obtain that
\begin{equation*}
\hat{f}_{a_1}\left(z^0_2,\ldots,z^0_N\right)=\hat{f}\left(a_1,z^0_2,\ldots,z^0_N\right),\qquad
a_1\in A_1\cap A^{\ast}_1.
\end{equation*}
This, combined with (\ref{eq8.2.4}), implies that
\begin{equation}\label{eq8.2.5}
\left\vert \hat{f} \left(\cdot,z^0_2,\ldots,z^0_N\right)\right\vert_{A_1\cap A^{\ast}_1}\leq\vert f\vert_A^{1-\delta}\vert f\vert_X^{\delta}
.
\end{equation}
On the other hand,
\begin{equation}\label{eq8.2.6}
\left\vert \hat{f} \left(\cdot,z^0_2,\ldots,z^0_N\right)\right\vert_{D_{1,\delta}}
\leq\vert \hat{f}\vert_{\widehat{X}}\leq\vert f\vert_X,
\end{equation}
where the latter estimate follows from Step 1.

Applying Theorem \ref{two-constant} to the function
 $\left.\log{\left\vert \hat{f} \left(\cdot,z^0_2,\ldots,z^0_N\right)\right\vert}\right|_{D_{1,\delta}},$ and taking (\ref{eq8.2.5}) and
(\ref{eq8.2.6}) into account, we obtain
\begin{eqnarray*}
\vert\hat{f}(z^0)\vert&\leq&
\left\vert \hat{f} \left(\cdot,z^0_2,\ldots,z^0_N\right)\right\vert_{A_1\cap
A^{\ast}_1}^{1- \widetilde{\omega}(z^0_1,A_1\cap A^{\ast}_1, D_{1,\delta})}
\left\vert \hat{f}
\left(\cdot,z^0_2,\ldots,z^0_N\right)\right\vert_{D_{1,\delta}}^{ \widetilde{\omega}(z^0_1,A_1\cap A^{\ast}_1,
 D_{1,\delta})}\\
&=&  \vert f\vert_A^{1-\omega(z^0)}\vert f\vert_X^{\omega(z^0)},  \end{eqnarray*}
where the equality follows from (\ref{eq8.2.3}) and the  identity
$  \widetilde{\omega}(z^0_1,A_1\cap A^{\ast}_1, D_{1,\delta})=\frac{ \widetilde{\omega}(z^0_1,A_1\cap A^{\ast}_1,
D_{1})}{1-\delta}$ (by Part 4) of Proposition \ref{prop3.5}). Hence estimate (\ref{eq8.2.2}) for the point
$z^0$ is proved. Since $z_0$ is an arbitrary point in $\widehat{X},$
(\ref{eq8.2.2}) follows, and the proof of the estimate in Theorem A is thereby finished.
\hfill $\square$

\smallskip

Combining the result of Subsections 7.1 and 7.2, Theorem A follows. \hfill
$\square$

\smallskip

Finally, we conclude the article by some remarks and open questions.


\smallskip

1. Recent development in the theory of separately analytic mappings is characterized by cross theorems with pluripolar
singularities and  boundary cross theorems.    The most general results are contained in some articles of Jarnicki and Pflug (see
\cite{jp2,jp3,jp4})
and in  recent works of  Pflug and the author  (see \cite{pn1,pn2}).
The question naturally arises whether one can generalize these results in
the context of mappings defined on complex manifolds with values in a
complex analytic space.  We postpone this issue to an ongoing work.

\smallskip

2. Is Theorem A in the case where $Z:=\C$ optimal? In other words, is the
open set $\widehat{X}$ always the envelope of holomorphy for separately holomorphic
functions defined on $X$?

\end{document}